\documentclass[a4paper,11pt,leqno]{article}
\usepackage{latexsym}
\usepackage{amsmath}
\usepackage{graphicx}
\usepackage{delarray}
\makeatletter
\@addtoreset{equation}{section}
\makeatother 
\usepackage{theorem}

\newtheorem{rem}{Remark}[section]
\newtheorem{cor}{Corollary}[section]
\newtheorem{thm}{Theorem}[section]

\newcommand{\qed}{\hfill $\Box$ \\ }
\theorembodyfont{\rmfamily}

\def\sgn{\mathop{\rm sgn}}

\title{Goodness-of-Fit Tests for Symmetric Stable Distributions 
-- Empirical Characteristic Function Approach}
\author{Muneya MATSUI${}^{\dag}$ and Akimichi TAKEMURA${}^{\ddag 
}$ \\ \\
$\dag$ Graduate School of Economics, University of Tokyo
 \\
$\ddag$ Graduate School of Information Science and Technology,
University of Tokyo
}
\date{October, 2005}

\begin{document}
\maketitle

\begin{abstract}
  We consider goodness-of-fit tests of symmetric stable distributions
  based on weighted integrals of the squared distance between the
  empirical characteristic function of the standardized data and the
  characteristic function of the standard symmetric stable
  distribution with the characteristic exponent $\alpha$ estimated
  from the data.  We treat $\alpha$ as an unknown parameter, but for
  theoretical simplicity we also consider the case that $\alpha$ is
  fixed.  For estimation of parameters and the standardization of data
  we use maximum likelihood estimator (MLE) and an equivariant
  integrated squared error estimator (EISE) which minimizes the
  weighted integral. We derive the asymptotic covariance function of
  the characteristic function process with parameters estimated by MLE
  and EISE.  
  For the case of MLE, the eigenvalues of the covariance function are
  numerically evaluated and asymptotic distribution of the test
  statistic is obtained using complex integration.
  Simulation studies show that the asymptotic distribution of the test
  statistics is very accurate.  
  We also present a formula of the asymptotic covariance function of
  the characteristic function process with parameters estimated by an
 efficient estimator for general distributions.
\end{abstract}

\section{Introduction.}
The family of stable distributions is one of the most important classes of
distributions in probability theory. The general central limit theorem
asserts that if a suitably normalized sum of independently and identically
distributed (i.i.d.) random variables has a limit distribution, only
possible limits are the stable distributions (Chapter 6 of Feller (1971)).
Concerning statistical inference, because of their attractive properties
such as heavy tails, many models based on stable distributions have been
considered in both social and natural sciences (Samorodnitsky and Taqqu
(1994), Uchaikin and Zolotarev (1999), Rachev and Mittnik (2000)).
Therefore it is important to consider goodness-of-fit tests of stable
distributions.  However few researches on goodness-of-fit tests of 
stable distributions have been conducted due to the difficulty in expressing
their density functions explicitly.  The purpose of this paper is to
propose goodness-of-fit tests based on the empirical characteristic
function, since the characteristic functions of stable distributions are
explicitly given. For past researches on goodness-of-fit tests of
heavy-tailed distributions using empirical characteristic function
approach, see G\"{u}rtler and\ Henze (2000) and Matsui and Takemura
(2005).
Both papers treat Cauchy ($\alpha=1$) distribution which
is one of the stable distributions.

Let $f(x;\mu,\sigma,\alpha)$ denote the symmetric stable
density with the characteristic function 
$$
\Phi(t)=
\exp(i\mu t-|\sigma t|^{\alpha}),
$$
where the parameter space is
$$
\Omega=\{-\infty <\mu<\infty,\ \sigma>0, \ 0<\alpha\leq 2 \}.
$$
Here $\alpha$ is the characteristic exponent, $\mu$ is the location
parameter and $\sigma$ is the scale
parameter. 
For the standard case $(\mu,\sigma)=(0,1)$ we simply write
the characteristic function as $\Phi(t;\alpha)=\exp(-|t|^\alpha)$ and the
density function as $f(x;\alpha)$. 
In this parameterization stable
distributions form a location-scale family for each value of $\alpha$, i.e.,
\[
f(x;\mu,\sigma,\alpha)=\frac{1}{\sigma}f(\frac{x-\mu}{\sigma};\alpha).
\]
In order to cope with more
general situation or for notational convenience we also write the parameters as 
\[
\theta=(\theta_1,\theta_2,\theta_3)=(\mu,\sigma,\alpha)  
\]
and write corresponding density, distribution or characteristic function
as 
\[
 f(x;\mu,\sigma,\alpha)=f(x;\theta),\quad
 F(x;\mu,\sigma,\alpha)=F(x;\theta),\quad \Phi(x;\mu,\sigma,\alpha)=\Phi(x;\theta).
\]
Here we note that $\theta$ is a vector.

In this paper we often differentiate functions of the parameter $\theta$
and the data $x$ with respect to $x,\mu,\sigma$ and $\alpha$.  
Since we will consider affine invariant (location-scale invariant) 
tests, it is often sufficient to evaluate the
derivatives at the standard case $(\mu,\sigma)=(0,1)$.  
For example we use the notation
\[
f_\mu(x;\alpha) = 
\frac{\partial}{\partial \mu}
f(x;\mu,\sigma,\alpha)_{|(\mu,\sigma)=(0,1)}\quad \mbox{or}\quad
f'(x;\alpha) = \frac{\partial}{\partial x} f(x;\mu,\sigma,\alpha)_{|(\mu,\sigma)=(0,1)} .
\]
Concerning the characteristic function we also use $\nabla_\theta
\Phi(t;\theta)=(\Phi_\mu(t;\theta),\Phi_\sigma(t;\theta),\Phi_\alpha(t;\theta))$
where, for example,
\[
 \Phi_\mu(t;\theta)=\frac{\partial}{\partial \mu}\Phi(t;\mu,\sigma,\alpha).
\]
For standard case $(\mu,\sigma)=(0,1)$ we write 
\[
 \Phi_\mu(t;\alpha)=\frac{\partial}{\partial \mu} \Phi(t;\mu,\sigma,\alpha)_{|(\mu,\sigma)=(0,1)}.
\]

Given a random sample $x_1,\ldots,x_n$ from an
unknown distribution $F$, we want to test the null hypothesis $H_1$ that $F$
belongs to the family of stable distributions $f(x;\mu,\sigma,\alpha)$
and the null hypothesis $H_2$ that $F$
belongs to the family of stable distributions $f(x;\mu,\sigma,\alpha)$
with $\alpha=\alpha_0$ fixed. Note that $H_1 \supset H_2$.  
Here we explain our proposed procedure for testing $H_1$, because for
$H_2$ we can simply replace $\hat\alpha$ by $\alpha_0$.

As remarked above stable distributions form a location scale family
and we consider affine invariant tests.  The proposed tests are based
on the difference between the empirical characteristic function
\begin{equation}
  \Phi_n(t) = \Phi_n(t;\hat\mu, \hat\sigma) 
            = \frac{1}{n}\sum \limits_{j=1}^{n}  \exp(ity_{j}), \qquad
  y_j = \frac{x_j - \hat\mu}{\hat\sigma}, 
 \label{eq:c.f.}
\end{equation}
of the standardized data $y_j$ and the characteristic
function with $\alpha$ estimated from the data
$$
\Phi(t)=\Phi(t;\hat\alpha)= e^{-{|t|}^{\hat\alpha}}.
$$
Here $\hat\mu=\hat{\mu}_n=\hat{\mu}_n(x_1,\ldots,x_n)$,
$\hat\sigma=\hat{\sigma}_n=\hat{\sigma}_n(x_1,\ldots,x_n)$
and $\hat\alpha=\hat{\alpha}_n=\hat{\alpha}_n(x_1,\ldots,x_n)$ 
are affine equivariant 
estimators  of $\mu$, $\sigma$, $\alpha$ satisfying 
\begin{eqnarray*}
\hat\mu_n(a+ b x_1, \ldots, a+ b x_n) 
  & = &  a + b \hat\mu_n(x_1,\ldots,x_n),\\
\hat\sigma_n(a+ b x_1, \ldots, a+ b x_n) 
  & = &   b \hat\sigma_n(x_1,\ldots,x_n),\\
\hat\alpha_n(a+ b x_1, \ldots, a+ b x_n) 
  & = &  \hat\alpha_n(x_1,\ldots,x_n),
\end{eqnarray*}
for all $-\infty < a <\infty$ and $b > 0$.

As equivariant estimators we 
consider maximum likelihood 
estimator (MLE) and an equivariant integrated squared error
estimator (EISE) defined in (\ref{eq:definition-EISE}) below. 
The reason for considering MLE is its asymptotic efficiency and the
reason for EISE is that its definition is similar to our proposed test
statistic.

Following G\"{u}rtler and Henze (2000) and 
Matsui and Takemura (2005) we propose the following
test statistic
\begin{equation}
  D_{n,\kappa}:=n\int_{-\infty}^{\infty}\big|\Phi_n(t)-e^{-{|t|}^{\hat\alpha}}\big|^2w(t)dt,
  \qquad w(t)=e^{-\kappa{|t|}}, \ \kappa >0. 
\label{eq:teststa}
\end{equation}
$D_{n,\kappa}$ is the weighted $L^2$-distance between $\Phi_n(t)$ and 
the characteristic function 
$e^{-{|t|}^{\hat\alpha}}$ of 
$f(x;\hat\alpha)$ with respect to the weight function 
$w(t)=e^{-\kappa{|t|}}, \ \kappa >0$.  This weight
function is chosen for convenience, so that we can evaluate
the asymptotic covariance function of the empirical characteristic
function process under $H_1$.

The test statistic $D_{n,\kappa}$ has an alternative representation, 
which is useful for obtaining its asymptotic distribution.
\begin{equation}
\label{eq:dkappa}
D_{n,\kappa}
= \int_{-\infty}^\infty {|\hat{Z}_n(t)|}^2\hat{\sigma}_n w(\hat{\sigma}_n t) dt,
\end{equation}
where
\begin{equation}
\hat{Z}_n(t) = \frac{1}{\sqrt{n}}\sum_{j=1}^{n} 
    \left \{\cos(tx_j)+i\sin(tx_j)-e^{-|\hat{\sigma}_nt|^{\hat{\alpha}_n}} 
\bigl( \cos(t \hat{\mu}_n) + i\sin(t \hat{\mu}_n) \bigl) \right \}.
\label{eq:ecp}
\end{equation}
$\hat{Z}_n(t)$ corresponds to  the empirical characteristic function process. 
 
Our test statistic $D_{n,\kappa}$ is a quadratic form of the empirical
characteristic function process.  Although we derive an explicit form
of the asymptotic covariance function of the empirical characteristic
function process, it is not trivial to derive the asymptotic
distribution of $D_{n,\kappa}$ under $H_1$ and $H_2$ from the
covariance function, especially when the parameters are estimated. 
See chapter 7 of Durbin (1973a) and Durbin (1973b) for tests based on
empirical distribution functions with estimated parameters and see
G\"{u}rtler and Henze (2000) and Matsui and Takemura (2005) for tests
based on empirical characteristic functions.
As in Matsui and Takemura (2005) we evaluate the asymptotic distribution
of $D_{n,\kappa}$  for the case MLE 
by numerically approximating the eigenvalues of the asymptotic covariance
function.  By numerically evaluating the asymptotic distribution
we can also check the convergence of the finite
sample distributions which we obtain by Monte Carlo simulations.

Concerning EISE, as shown below, the asymptotic covariance function of
the empirical characteristic function process is very complicated.  
Furthermore we found that Monte Carlo simulation involving EISE is
very time consuming.  Therefore in this paper we show theoretical
results on our proposed test statistic involving EISE and leave
numerical studies to our subsequent works.

This paper is organized as follows. In Section \ref{sec:estimator} we
first define and summarize properties of MLE and EISE.  Then in
Section \ref{theory-of-test-statstics}
we state theoretical results on asymptotic distribution of $D_{n,\kappa}$ under
$H_1$ in Theorem \ref{thm:2.3} and Theorem
\ref{thm:weak-convergence-of-test-statistics} and results 
under $H_2$ as the corollaries of these  theorems.
Numerical evaluations of asymptotic
critical values of $D_{n,\kappa}$ under $H_1$ and $H_2$ for MLE are
discussed in Section
\ref{sec:approximation-of-test-asymptotic-statistic}. Simulation
studies of MLE and corresponding test statistic $D_{n,\kappa}$ are
given in Section 4, including the study of finite sample power behavior 
in Section \ref{sec:power}.

\section{Main results}

\subsection{Estimators and their asymptotic properties}
\label{sec:estimator}
For our purposes we need asymptotic covariance matrices 
and ``asymptotically linear representations'' (AL representations) of the 
estimators.
We describe asymptotic properties of maximum likelihood estimator (MLE) 
following DuMouchel (1973). 
We also define an equivariant 
integrated squared error estimator (EISE) and give
asymptotic properties of EISE.
For MLE explicit expressions of the asymptotic covariance matrix 
and AL representations are given in the Cauchy case $\alpha=1$.

As shown in DuMouchel (1973), MLE is asymptotically normal
and asymptotically  efficient. The likelihood equation is given by
\begin{eqnarray}
\label{eq:likelihood-equation-1}
\frac{\partial L}{\partial \mu} = 0 
  &\Leftrightarrow & 
    \sum_{j=1}^{n} \frac{1}{2\pi f(\frac{x_j-\mu}{\sigma};\alpha)}
\int_{-\infty}^{\infty} 
e^{-it\left(\frac{x_j-\mu}{\sigma}\right)} \Phi_\mu(t;\alpha)
dt= 0,  \\
\label{eq:likelihood-equation-2}
\frac{\partial L}{\partial \sigma} = 0 
&\Leftrightarrow & 
\sum_{j=1}^{n} \frac{1}{2\pi f(\frac{x_j-\mu}{\sigma};\alpha)}
\int_{-\infty}^{\infty}  
e^{-it\left(\frac{x_j-\mu}{\sigma}\right)} \Phi_\sigma(t;\alpha)
dt= 0,\\
  \label{eq:likelihood-equation-3}
  \frac{\partial L}{\partial \alpha} = 0 
  &\Leftrightarrow & 
  \sum_{j=1}^{n} \frac{1}{2\pi  f(\frac{x_j-\mu}{\sigma};\alpha)}
  \int_{-\infty}^{\infty} 
e^{-it\left(\frac{x_j-\mu}{\sigma}\right)}
\Phi_\alpha(t;\alpha)
dt= 0,  
\end{eqnarray}
where 
\begin{equation*}
\left(\Phi_\mu(t;\alpha),\Phi_\sigma(t;\alpha),\Phi_\alpha(t;\alpha)\right)=
  \left(ite^{-|t|^\alpha},-e^{-|t|^\alpha}|t|^{\alpha}\alpha,
       -e^{-|t|^\alpha}|t|^{\alpha}\log|t| \right).
\end{equation*}

EISE is an affine equivariant version of the ISE (integrated squared
error) estimator proposed by Paulson \textit{et al.} (1975). The
original ISE estimator of Paulson \textit{et al.} (1975) is not
equivariant.  Robustness and efficiency of ISE estimators of location
and scale parameters are discussed in Thornton and Paulson
(1977) for the normal case and in Besbeas and Morgan (2001) for the Cauchy case.
EISE is based on the standardized empirical characteristic function.
Let 
\[
\Phi_n(t; \mu,\sigma) = \frac{1}{n}\sum_{j=1}^n
\exp\Big(it\frac{x_j-\mu}{\sigma}\Big), 
\]
which is the same as (\ref{eq:c.f.}) with
$\hat{\mu}_{n}$ and $\hat{\sigma}_n$ replaced by $\mu$ and $\sigma$.
Write
\begin{equation}
\label{eq:eise}
Q(\mu, \sigma, \alpha)
=\int_{-\infty}^{\infty}{|\Phi_n(t; \mu, \sigma)-e^{-{|t|}^\alpha
}|}^2  w(t)dt,
\end{equation}
where we use the following weight function 
\begin{equation}
\label{eq:ise-equivaiant-weight}
   w(t)= \exp(-\nu{|t|}^{\bar{\alpha}}),\quad \nu>0.
\end{equation}
Here we call $\bar{\alpha}$ weighting index and $\nu$ weighting constant.
EISE $(\hat\mu_n, \hat\sigma_n, \hat\alpha_n)$ is defined to be the
minimizer of 
$Q(\mu, \sigma, \alpha)$:
\begin{equation}
\label{eq:definition-EISE}
Q(\hat\mu_n, \hat\sigma_n, \hat\alpha_n)
 = \min_{\mu, \sigma, \alpha} Q(\mu,\sigma,\alpha).
\end{equation}
It is easy to see that  EISE is affine equivariant by definition.
Note that the weighting constant $\kappa$ in
the test statistic (\ref{eq:teststa}) and the weighting constant $\nu$ 
in (\ref{eq:ise-equivaiant-weight}) for
EISE may be different.  In our theoretical results on EISE 
we can treat more general weighting functions, i.e., $w(t) \ge 0$ is an
arbitrary even function.
However for performing goodness-of-fit tests, it seems natural to set
$\alpha_0=\bar\alpha$ and $\nu=\kappa$. 
The integral $Q(\theta)$ can be calculated as
$$
Q(\theta)=\int_{-\infty}^{\infty}\biggl\{
\frac{1}{n^2}\sum_{j,k}^{n}\cos\big(t(x_j-x_k)/\sigma\big)
-\frac{2}{n}\sum_{j=1}^{n}\cos\big(t(x_j-\mu)/\sigma
\big)e^{-|t|^{\alpha}}+e^{-2|t|^{\alpha}}
 \biggr\} w(t)dt.
$$
The estimators satisfy the following
estimating equations $0=\partial Q/\partial \mu = \partial Q/\partial\sigma =\partial
Q/\partial \alpha \Leftrightarrow 
0=Q_\mu(\theta)=Q_\sigma(\theta)=Q_\alpha(\theta)$. 
\begin{eqnarray}
&&Q_\mu(\theta)
= 
-\int_{-\infty}^{\infty} 
\biggr \{ \frac{1}{n}\sum_{j=1}^{n}\sin\big(t(x_j-\mu)/\sigma\big)
t e^{-|t|^{\alpha}} \biggl \} w(t)dt, \label{eq:eisemu}\\
&&Q_\sigma(\theta)
=
\frac{1}{2}
\int_{-\infty}^{\infty}
\biggr [ \biggr \{\frac{1}{n^2}\sum_{j,k=1}^n \cos\big(t(x_j-x_k)\big)
- \frac{2}{n} \sum_{j=1}^{n}\cos\big(t(x_j-\mu)\big)e^{-|\sigma t|^{\alpha}}
+ e^{-2|\sigma t|^{\alpha}} \biggl \} \label{eq:eisesig}\\
&&\hspace{2.0cm} \times
\{w'(\sigma |t|)\sigma|t|+ w(\sigma t)\} +
2\biggr \{\frac{1}{n} \sum_{j=1}^{n}
\cos\big(t(x_j-\mu)\big)
-e^{-|\sigma t|^{\alpha}} \biggl \}
\alpha|\sigma t|^{\alpha} e^{-|\sigma t|^{\alpha}}
w(\sigma t)\biggl ] dt, \nonumber \\
&&Q_\alpha(\theta)
= 
\int_{-\infty}^{\infty}
\biggr \{\frac{1}{n} \sum_{j=1}^{n}\cos\big(t(x_j-\mu)/\sigma\big)
-e^{-|t|^{\alpha}} \biggl \}
e^{-|t|^{\alpha}}|t|^{\alpha}\log|t| w(t)dt,\label{eq:eiseal}   
\end{eqnarray}
where $w'(x)=dw(x)/dx$.  Note that in case of $Q_\sigma(\theta)$
differentiation was done after the transformation $t \rightarrow \sigma t$.

In the rest of this paper we use the following notations.
$\stackrel{D}{\longrightarrow}$
means weak convergence of random variables or stochastic processes,
$\stackrel{P}{\longrightarrow}$ means convergence in probability.
 
The asymptotically linear representations (AL representation) give an
method of approximating asymptotic behavior of the estimator by sum of
functions of i.i.d.\ random samples. For the standard symmetric stable case
$f(x;\alpha)$ we need following three expressions,
\begin{eqnarray*}
\quad \sqrt{n}\hat{\mu}_n
&=& \frac{1}{\sqrt{n}}\sum_{j=1}^{n}l_1(X_j) 
 + r_{1n} , \\
\quad \sqrt{n}(\hat{\sigma}_n-1)
&=& \frac{1}{\sqrt{n}}\sum_{j=1}^{n}l_2(X_j) 
 + r_{2n}, \\
\sqrt{n}(\hat{\alpha}_n-\alpha) 
&=& \frac{1}{\sqrt{n}}\sum_{j=1}^{n}l_3(X_j) 
  + r_{3n},\qquad  r_{1n}, r_{2n}, r_{3n} \stackrel{P}{\longrightarrow} 0. \label{baha} 
\end{eqnarray*}

For the case of MLE, AL representations are given in terms of the score
functions ((\ref{eq:likelihood-equation-1}-\ref{eq:likelihood-equation-3}))
and the Fisher information matrix.  The proof is standard and omitted.

\begin{thm}
\label{thm:2.1mle}
Let $I(\theta)$ denote the Fisher 
information matrix
\[
I(\theta)=
\left(
\begin{array}{ccc}
 I_{11}     & 0       & 0       \\
  0         & I_{22} & I_{23}  \\
  0         & I_{32}  & I_{33}
\end{array}
\right), \qquad
I_{ij}(\theta)=
- \mbox{E}_\theta \left[\frac{\partial^2 \log f(x,\theta)}
{\partial \theta_i \partial \theta_j}\right].
\]
The AL representations $l_{\theta}(x)=(l_\mu(x),l_\sigma(x),l_\alpha(x))$
at the standard case $(\mu,\sigma)=(0,1)$
are given by $l_{\theta}(x)=I^{-1}(\theta) h_{\theta}(x)$, where $h_\theta(x)$
are
\begin{eqnarray} 
h_{\mu}(x)&=&
\frac{1}{2\pi f(x;\alpha)}
\int_{-\infty}^{\infty} it e^{-itx }e^{-{|t|}^\alpha}
dt, \label{mle:Bahadurmu}\\
h_{\sigma}(x)&=&
-\frac{\alpha}{2\pi f(x;\alpha)}
\int_{-\infty}^{\infty} e^{-itx} e^{-{|t|}^\alpha}
{|t|}^{\alpha} dt\label{mle:Bahadursigma},\\
h_{\alpha}(x)&=&
-\frac{1}{2\pi f(x;\alpha)}
\int_{-\infty}^{\infty} e^{-itx}e^{-{|t|}^\alpha}
{|t|}^{\alpha}\log|t| dt.\label{mle:Bahaduralpha} 
\end{eqnarray}
\end{thm}
Concerning EISE we can employ standard theory of $U$-statistics.
The proof is given in Appendix \ref{sec:a1}.
  
\begin{thm}\label{thm:2.1eise}
Define a $3\times 3$ symmetric matrix
\begin{equation}
\label{eq:matrix-A}
A=A(\alpha)=A(\theta)_{|(\mu,\sigma)=(0,1)}
\end{equation}
by
\begin{eqnarray*}
A_{12}&=& A_{12}=0, \\
A_{11}
&=& 
\int^{\infty}_{-\infty}e^{-2{|t|}^{\alpha}}t^2 w(t)dt, \\
A_{22}
&=& \alpha^2 \int^{\infty}_{-\infty}e^{-2{|t|}^{\alpha}}{|t|}^{2\alpha} w(t)dt, \\
A_{23}
&=& \alpha \int^{\infty}_{-\infty}e^{-2{|t|}^{\alpha}}{|t|}^{2\alpha}
\log|t| w(t)dt, \\
A_{33}
&=& \int^{\infty}_{-\infty}e^{-2{|t|}^{\alpha}}{|t|}^{2\alpha}
{(\log|t|)}^2 w(t)dt,
\end{eqnarray*}
and define $h_{\theta}(x)=(h_{\mu}(x),h_{\sigma}(x),h_{\alpha}(x))'$ by 
\begin{eqnarray}
h_{\mu}(x)&=&
\int_{-\infty}^{\infty} t\sin(tx)e^{-{|t|}^\alpha}
w(t)dt, \\
h_{\sigma}(x)&=&
-\alpha\int_{-\infty}^{\infty} (\cos(tx)-e^{-{|t|}^\alpha})e^{-{|t|}^\alpha}
{|t|}^{\alpha} w(t)dt, \\
h_{\alpha}(x)&=&
-\int_{-\infty}^{\infty}(\cos(tx)-e^{-{|t|}^\alpha})
e^{-{|t|}^\alpha}{|t|}^{\alpha}\log|t| w(t)dt. 
\end{eqnarray}
For EISE the AL representations $l_{\theta}(x)=(l_\mu(x),l_\sigma(x),l_\alpha(x))$ 
at the standard case $(\mu,\sigma)=(0,1)$ are given by
\begin{equation}
l_{\theta}(x)=A(\alpha)^{-1} h_{\theta}(x)
\end{equation}
and their asymptotic covariance matrix at the standard case is given by 
\[
A^{-1}\mbox{E}[h_\theta(X)h_\theta'(X)]A^{-1}{}'=
A^{-1}\left(
\begin{array}{ccc}
H_{\mu\mu}       &  0               & 0                 \\
0                & H_{\sigma\sigma} & H_{\sigma\alpha}  \\
0                & H_{\alpha\sigma} & H_{\alpha\alpha} \\
\end{array}
\right)A^{-1}{}'
=
\left(
\begin{array}{ccc}
J_{11} &  0        & 0       \\
0      & J_{22}    & J_{23}  \\
0      & J_{32}    & J_{33}  \\
\end{array}
\right),
\]
where each element of $H=\mbox{E}[h_\theta(X)h_\theta'(X)]$ is
\begin{eqnarray*}
H_{\mu\mu}
&=& 
\int_{-\infty}^{\infty}\int_{-\infty}^{\infty}
\left \{
\frac{1}{2}(e^{-{|s-t|}^\alpha}-e^{-{|s+t|}^\alpha})
\right \}
e^{-(|s|^\alpha+|t|^\alpha)}st\ w(s)w(t)dsdt, \\
H_{\sigma \sigma}
&=& 
\alpha^2\int_{-\infty}^{\infty}\int_{-\infty}^{\infty}
\left \{
\frac{1}{2}(e^{-{|s-t|}^\alpha}+e^{-{|s+t|}^\alpha})-e^{-(|s|^\alpha+|t|^\alpha)}
\right \}
e^{-(|s|^\alpha+|t|^\alpha)}{|st|}^\alpha w(s)w(t)dsdt, \\
H_{\sigma\alpha}
&=& 
\alpha\int_{-\infty}^{\infty}\int_{-\infty}^{\infty}
\left \{
\frac{1}{2}(e^{-{|s-t|}^\alpha}+e^{-{|s+t|}^\alpha})-e^{-(|s|^\alpha+|t|^\alpha)}
\right \}
e^{-(|s|^\alpha+|t|^\alpha)}{|st|}^\alpha \log|t|w(s)w(t)dsds, \\
H_{\alpha\alpha}
&=& 
\int_{-\infty}^{\infty}\int_{-\infty}^{\infty}
\left \{
\frac{1}{2}(e^{-{|s-t|}^\alpha}+e^{-{|s+t|}^\alpha})-e^{-(|s|^\alpha+|t|^\alpha)}
\right \}
e^{-(|s|^\alpha+|t|^\alpha)}{|st|}^\alpha \log|s|\log|t|w(s)w(t)dsdt  .
\end{eqnarray*}
\end{thm}

Note that the above AL representations and asymptotic matrices involve definite
integrals, which require numerical integration. 
But for some special cases like Cauchy ($\alpha =1$) 
we can calculate several integrals analytically.
Analytic expressions are useful for checking correctness of
numerical calculations concerning Theorem \ref{thm:2.1mle}.
We give the following Corollary for the case of $\alpha=1$ and MLE.

\begin{cor}\label{cor:2.1}
Let $\gamma \doteq 0.577216$ denote Euler constant.
In the Cauchy case ($\alpha=1$) and MLE, at $(\mu,\sigma)=(0,1)$, 
the AL representations 
are given as $l_\theta(x)=I^{-1} h_\theta(x)|_{\alpha=1},$ where
\begin{eqnarray*}
h_\mu(x)
&=&\frac{2x}{x^2+1} ,\qquad 
h_\sigma(x)
=\frac{x^2-1}{x^2+1},\\
h_\alpha(x)
&=&\frac{1-x^2}{x^2+1}
\left[ \frac{1}{2}
\log(x^2+1)-1+\gamma \right]
+\frac{2x}{x^2+1}
\arctan x,
\end{eqnarray*}
\[
I_{11}=I_{22}= \frac{1}{2},
\quad I_{23}=I_{32}=
\frac{1}{2}(1-\gamma-\log2),\quad
 I_{33}=\frac{1}{2}
\left\{\frac{\pi^2}{6}+{(\gamma+\log2-1)}^2\right\}.
\]
\end{cor}

A similar result is given in Section 6 of Matsui and Takemura (2006)
and the proof is omitted.

\subsection{Asymptotic theory of the proposed test statistics}
\label{theory-of-test-statstics}
 In this section theoretical results on asymptotics
 of the proposed test statistic $D_{n,\kappa}$ are obtained.
 {}From another expression of $D_{n,\kappa}$ 
 (\ref{eq:dkappa}) we derive weak convergence of $\hat{Z}_n(t)$ 
 and weak convergence of test statistic $D_{n,\kappa}$ in the following
 two theorems. These results correspond to those  of
 Cauchy case stated in Matsui and Takemura (2005) where parameter $\alpha=1$
 is fixed. As a special case we also describe the Cauchy case involving
 estimation of $\alpha$ in the corollary below. Furthermore a general
 formula of asymptotic covariance function of the empirical characteristic
 process with parameters estimated by an efficient estimator is given in
 the latter part of this section.
 As already remarked several times, we can assume without loss of generality that
 $X_1,\ldots,X_n$ is random sample from $f(x;\alpha)$ because of affine
 invariance of our tests.
 Following G\"{u}rtler and Henze (2000) we use 
 the Fr\'{e}chet space $C(\mathbf{R})$ of continuous functions on
 the real line $\mathbf{R}$ for considering the random processes.
 The metric of $C(\mathbf{R})$ is given by
\[
\rho(x,y)=\sum_{j=1}^{\infty}2^{-j}\frac{\rho_j(x,y)}{1+\rho_j(x,y)},
\]
where $\rho_j(x,y)=\max_{|t|\le j}|x(t)-y(t)|$. 

We first give the asymptotic covariance function of the empirical
characteristic function process with parameters estimated by MLE and EISE.
In the following theorems the elements of the inverse of the Fisher
information matrix $I(\theta)$  and the matrix $A$ in (\ref{eq:matrix-A})
are denoted with superscripts $I^{ij}$ and $A^{ij}$.

\begin{thm}
\label{thm:2.3}
  Let $X_1,\ldots,X_n$ be {\rm i.i.d.} $f(x;\alpha)$ random variables and
  let $\hat{Z_n}$ be defined in $( \ref{eq:ecp})$. Then $\hat{Z_n}$
  $\stackrel{D}{\longrightarrow} Z$ in $\mbox{C}(\mathbf{R})$, where
  $Z$ is a zero mean Gaussian process with covariance functions given 
  below.\\
{\rm MLE\ :} 
\vspace{-3mm}
\begin{eqnarray}
\label{eq:mle-covariance}
\Gamma(s,t) &=& e^{-{|t-s|}^\alpha}-e^{-({|t|}^\alpha+{|s|}^\alpha)} 
-\bigg \{ I^{11}st+ I^{22}{|st|}^{\alpha}\alpha^2  \\
& & +I^{23}{|st|}^{\alpha}\alpha\log|st| + I^{33}{|st|}^{\alpha}\log|s|\log|t|
\bigg\}e^{-({|t|}^\alpha+{|s|}^\alpha)} \nonumber
.\end{eqnarray}
{\rm EISE\ :}
\vspace{-3mm}
{\allowdisplaybreaks
\begin{align}
\Gamma(s,t)&=
e^{-{|t-s|}^\alpha}-e^{-({|t|}^\alpha+{|s|}^\alpha)} \nonumber \\
& \quad +
\{J_{11}st
+J_{22}\alpha^2{|st|}^\alpha
+J_{23}\alpha{|st|}^\alpha(\log|s|+\log|t|)
+J_{33}{|st|}^\alpha \log|t| \log |s| \}
e^{-({|t|}^\alpha+{|s|}^\alpha)} \nonumber\\
& \quad +
\left\{(B_\sigma A^{22}+B_\alpha A^{23})\alpha
  ({|t|}^\alpha+{|s|}^\alpha)+(B_\sigma A^{23}+B_\alpha A^{33}
)({|t|}^\alpha \log|t|
+{|s|}^\alpha \log|s|)\right\}
e^{-({|t|}^\alpha+{|s|}^\alpha)} \nonumber\\
& \quad -
A^{11}t e^{-{|t|}^\alpha}
\int_{-\infty}^{\infty}
e^{{|s-u|}^\alpha-{|u|}^\alpha}
u\> w(u)du[2]
\nonumber \\
& \quad -
\alpha(A^{22}\alpha+A^{23}\log|t|)
{|t|}^\alpha e^{-{|t|}^\alpha}
\int_{-\infty}^{\infty}
e^{{|s-u|}^\alpha-{|u|}^\alpha}
{|u|}^\alpha w(u)du[2] \nonumber\\
& \quad-
(A^{23}\alpha+A^{33}\log|t|){|t|}^\alpha e^{-{|t|}^\alpha}
\int_{-\infty}^{\infty}
e^{{|s-u|}^\alpha-{|u|}^\alpha}
{|u|}^\alpha \log |u| w(u)du[2]
\label{eq:eise-covariance}
\end{align}
}
where $[2]$ after a term means symmetrization with respect to $s$ and $t$, i.e.,
$g(s,t)[2]=g(s,t)+g(t,s)$,  and
\begin{equation*}
B_\alpha=\int_{-\infty}^{\infty}e^{-2{|u|}^\alpha}{|u|}^\alpha \log |u| w(u)du,
\qquad
B_\sigma=\alpha\int_{-\infty}^{\infty}e^{-2{|u|}^\alpha}{|u|}^\alpha w(u)du.
\end{equation*}
\end{thm}
As a corollary the asymptotic covariance function
for the case of fixed $\alpha$ is given as follows.

\begin{cor}
Under the same conditions of Theorem \ref{thm:2.3}, when $\alpha$ is fixed,
$\hat{Z}_n(t)$$\stackrel{D}{\longrightarrow} Z$ in $\mbox{C}(\mathbf{R})$, where
  $Z$ is a zero mean Gaussian process with covariance functions given 
  below.\\
{\rm MLE\ :} 
\begin{equation}
\Gamma(s,t) = e^{-{|t-s|}^\alpha}-e^{-({|t|}^\alpha+{|s|}^\alpha)} 
-( I^{11}st+ I^{22}{|st|}^{\alpha}\alpha^2 )e^{-({|t|}^\alpha+{|s|}^\alpha)} \nonumber
.
\end{equation}
{\rm EISE\ :}
\begin{equation*}
\begin{split}
\Gamma(s,t) &=
 e^{-{|t-s|}^\alpha}-e^{-({|t|}^\alpha+{|s|}^\alpha)} 
+
\{J_{11}st+J_{22}\alpha^2{|st|}^\alpha+A^{22}B_\sigma\alpha(|t|^\alpha+|s|^\alpha)
    \}e^{-({|t|}^\alpha+{|s|}^\alpha)} \\
& \quad-A^{11}te^{-|t|^\alpha}\int_{-\infty}^{\infty}
  e^{-|s-u|^\alpha-|u|^\alpha}
u w(u)du[2]
-A^{22}\alpha^2 
|t|^\alpha e^{-|t|^\alpha}  \int_{-\infty}^{\infty}
e^{-|s-u|^\alpha-|u|^\alpha}
|u|^\alpha w(u)du[2].
\end{split}
\end{equation*}
\end{cor}

The following corollary gives the asymptotic covariance function when the
true distribution is Cauchy $C(\mu,\sigma)$, but the characteristic exponent $\alpha$ is
estimated by MLE.

\begin{cor}
  Let $X_1,\ldots,X_n$ be {\rm i.i.d.} $C(0,1)$ random variables and
  let $\hat{Z_n}$ be defined in $( \ref{eq:ecp})$, where parameters
  are estimated by MLE. Then $\hat{Z_n}$
  $\stackrel{D}{\longrightarrow} Z$ in $\mbox{C}(\mathbf{R})$, where
  $Z$ is a zero mean Gaussian process with covariance function given 
  below. 
\begin{eqnarray}
&& \Gamma_c (s,t)= e^{-|t-s|}-\{ 1+2(st+|st|) \} e^{-|s|-|t|} \\  
&& \hspace{1.5cm}-\frac{12}{\pi^2}\left\{\log|s|+(\gamma+\log 2-1)\right\} 
\left\{ \log|t|+(\gamma+\log 2-1) \right\}|st|e^{-|s|-|t|} .
\nonumber
\end{eqnarray}
\end{cor}

Furthermore for efficient estimations including MLE, we can derive
more general result after some formulations. We assume parameter space
$\Theta$ is $p$-dimensional.
First, we define an efficient estimator $\hat\theta_n$ of $\theta_0$ as 
such that
\begin{eqnarray}
\label{eq:AL-efficient}
\sqrt{n}(\hat{\theta}_n-\theta_0) &=& \frac{1}{\sqrt{n}}\sum_{j=1}^n
 I^{-1}(\theta_0)\frac{\partial f(X_j;\theta_0)}{\partial
 \theta}+\epsilon_n \nonumber \\
&:=& \frac{1}{\sqrt{n}}\sum_{j=1}^n
 l_E(X_j;\theta_0)+\epsilon_n, \label{eq:AL-efficient}
\end{eqnarray}
where $\epsilon_n \stackrel{P}{\longrightarrow} 0$. This definition
coincides with formula (20) of Durbin (1973b), though his definition
is given for nuisance parameters.
He also gives conditions that the estimator $\hat{\theta}_n$ satisfies formula
(\ref{eq:AL-efficient}) under more general arguments including local
alternatives (see the condition (A3) of Durbin (1973b)). 
Note that Durbin (1973b) considers estimation only for nuisance
parameters $\hat\theta_n$ including local alternatives and the
parameters of interest (the null hypothesis) are not estimated. However in our case
the whole parameters are estimated without nuisance parameters because
the null hypothesis of $H_1$ is the whole parameter space.

Second, we modify the conditions (iv) of Cs\"{o}rg\H{o} (1983)
as \\
$(iv)^*$  
\begin{itshape}
$l_E(x;\theta)$ in (\ref{eq:AL-efficient})
 is a $p$-dimensional Borel measurable function, 
 $E[l_E(X_1;\theta_0)]=(0,\ldots,0)$, and 
 $I^{-1}(\theta_0)=E[l_E(X_1;\theta_0)l'_E(X_1;\theta_0)]$ is
 finite and positive definite. 
\end{itshape}

\begin{thm}
\label{thm:efficient}
  Let $X_1,\ldots,X_n$ be {\rm i.i.d.}\ $F(x;\theta_0)$ random variables
 and let $k(x,t)=\cos(tx)+i\sin(tx)$. Consider 
 the kernel transformed empirical characteristic process 
\begin{equation}
\hat{Z_n}=
\int k(x,t) d\left\{\sqrt{n}\bigl(F_n(x)
   - F(x;\hat{\theta}_n)\bigl)\right\}, 
\end{equation}
 where $F_n(x)$  denotes the empirical distribution function.
 Then $\hat{Z_n}$
  $\stackrel{D}{\longrightarrow} Z$ in $\mbox{C}(\mathbf{R})$
 under the conditions $(i)^*$, $(ii)^*$, $(v)$ and $(vi)$ of
 Cs\"{o}rg\H{o} (1983) where $l(\cdot;\theta_0)$ is replaced by
 $l_E(\cdot;\theta_0)$ and $(iv)^*$.  
  Here $Z$ is a zero mean Gaussian process with covariance function
\begin{equation}
\Gamma(s,t) = \Phi(s-t;\theta_0)-\Phi(s;\theta_0)\overline{\Phi(t;\theta_0)}
            -\nabla_\theta\Phi(s;\theta_0)'I^{-1}(\theta)
              \overline{\nabla_\theta\Phi(t;\theta_0)}, 
\end{equation}
where $\Phi(t;\theta_0)$ is the characteristic function and
 $\nabla_\theta\Phi(t;\theta_0)$ is the derivative of $\Phi(t;\theta_0)$
 with respect to parameter vector $\theta$.
\end{thm}
Theorem \ref{thm:efficient} is interpreted as the Fourier kernel
transformed version of Theorem 2 of Durbin (1973b) with nuisance
parameters corresponds to the estimated null hypothesis. In our subsequent works
we will consider extension of Theorem \ref{thm:efficient} to the case of 
local alternatives.

Finally we state the following theorem concerning the
 weak convergence of
$D_{n,\kappa}$.

\begin{thm}
\label{thm:weak-convergence-of-test-statistics}
Under the conditions of Theorem 2.1

\[
D_{n,\kappa}=\int_{-\infty}^\infty \hat{Z}_n(t)^2\hat{\sigma}_n
      e^{-\hat{\sigma}_n\kappa {|\hat{\sigma}_nt|}^{\hat{\alpha}_n}} dt\stackrel{D}{\longrightarrow}
      D_\kappa:=\int_{-\infty}^\infty Z(t)^2e^{-\kappa{|t|}^\alpha} dt.
\]
\end{thm} 

The proofs of the above theorems are given in Appendix
\ref{sec:proof-of-section2.2}.  We omit the proof of Theorem
\ref{thm:weak-convergence-of-test-statistics} since after obtaining
Theorem \ref{thm:2.3} the proof is the same as that of Theorem 2.2 of
G\"{u}rtler and Henze (2000).

\section{Approximation of the asymptotic critical values of the proposed
  test statistics}
\label{sec:approximation-of-test-asymptotic-statistic}
In this section we investigate the distribution of $D_\kappa$ for
MLE. We briefly explain how to obtain the characteristic
function of $D_\kappa$. 
Detailed treatments of this approach in statistical
applications are given in Tanaka (1996) or Anderson and Darling (1952).
Since the characteristic function of $D_\kappa$ contains infinite
product of functions of eigenvalues which can not be evaluated analytically,
we approximate eigenvalues by theory of homogeneous integral equations
of the second kind and the associated Fredholm determinant.
Then utilizing complex integration, we invert the  characteristic function
and obtain series representation of the distribution of
$D_\kappa$. Detailed theoretical argument of inversion process is given
in Slepian (1957). Actual computational approximations are given in the
next section. At the end of this section we transform our kernels $\Gamma(s,t)$ on
$\mathbf{R}^2$ to kernels $K(s,t)$ on $[-1,1]^2$ for convenience in
numerical computation.

In this paper we omit consideration of $D_\kappa$ of EISE since kernels
$\Gamma(s,t)$ have definite integrals and we need many numerical
approximations. 
On the other hand for the case of MLE we can utilize past
researches on Fisher information in DuMouchel (1975), Nolan (2001) and
Matsui and Takemura (2006) to confirm the accuracy of our
computation.
For $\kappa > 1$ we can use the following standard form of Mercer's
theorem.

\begin{figure}
\begin{center}
\begin{minipage}{.45\linewidth}
\includegraphics[width=\linewidth]{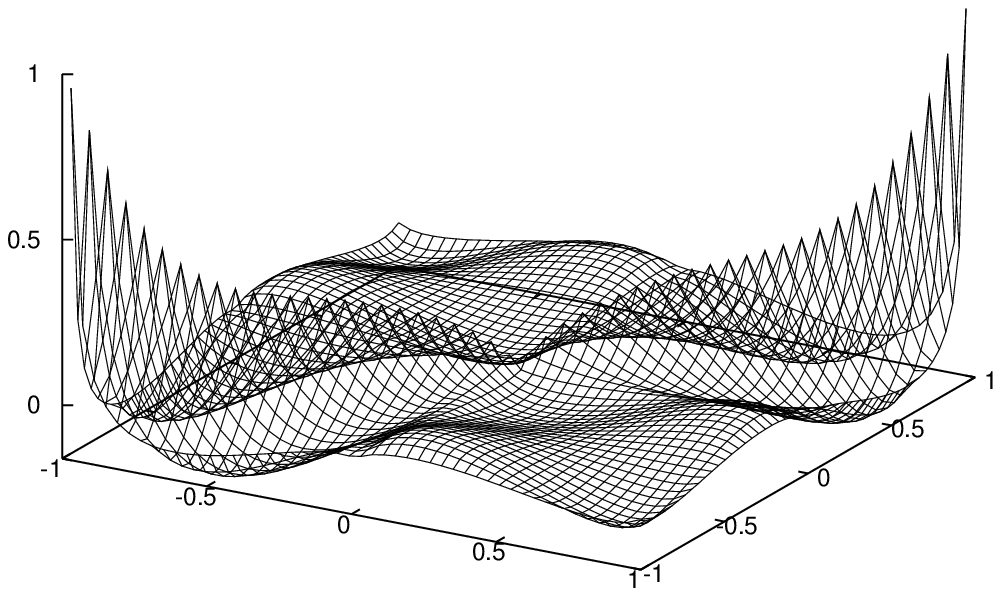}
\caption{MLE-H1 ($\alpha=1.0$, $\kappa=1.0$)}
\label{fig:cov-mle-a1.0-k1.0}
\end{minipage}
\begin{minipage}{.45\linewidth}
\includegraphics[width=\linewidth]{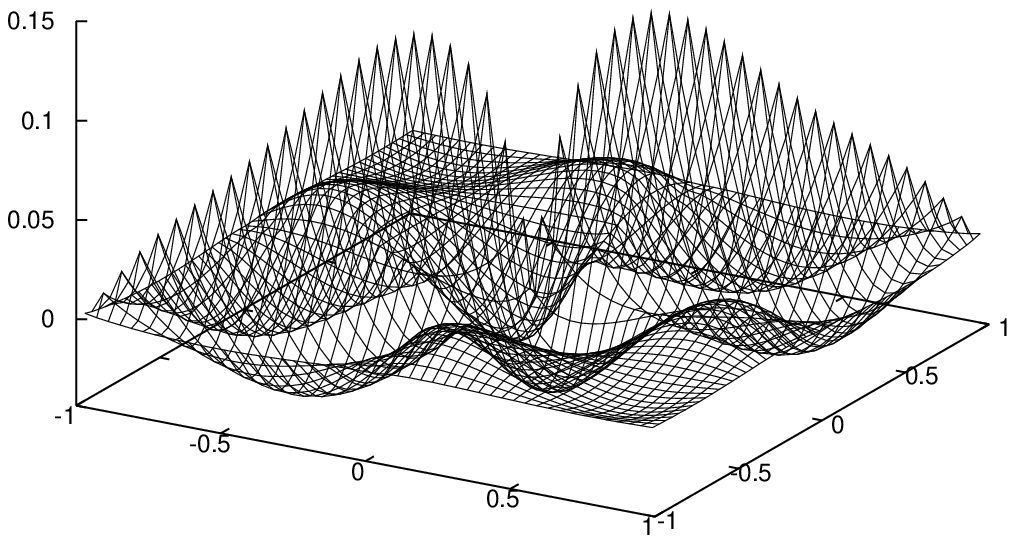}
\caption{MLE-H1 ($\alpha=1.0$, $\kappa=2.5$)}
\label{fig:cov-mle-a1.0-k2.5}
\end{minipage}
\end{center}
\end{figure}

\begin{figure}
\begin{center}
\begin{minipage}{.45\linewidth}
\includegraphics[width=\linewidth]{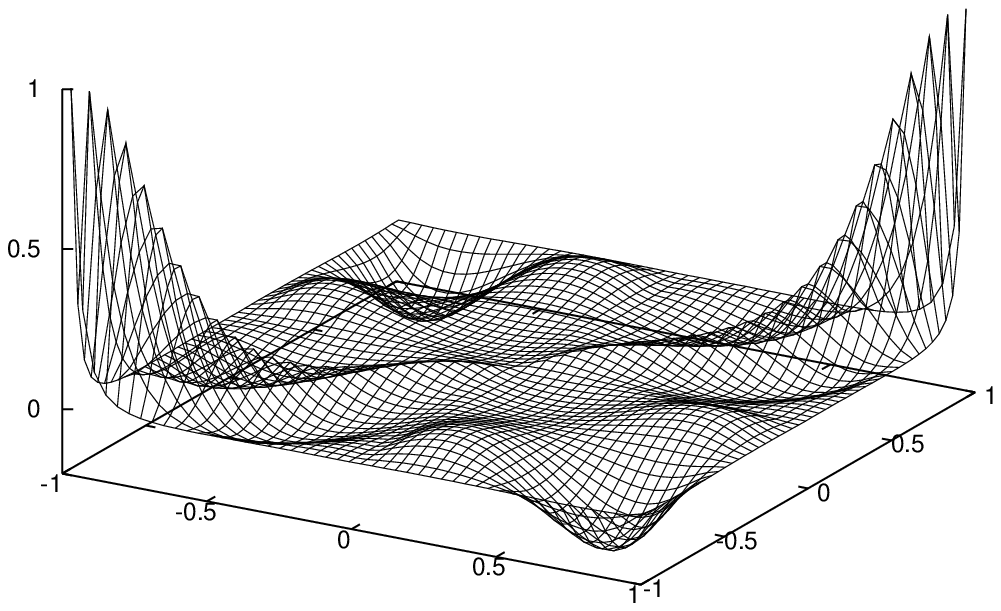}
\caption{MLE-H1 ($\alpha=1.5$, $\kappa=1.0$)}
\label{fig:cov-mle-a1.5-k1.0}
\end{minipage}
\begin{minipage}{.45\linewidth}
\includegraphics[width=\linewidth]{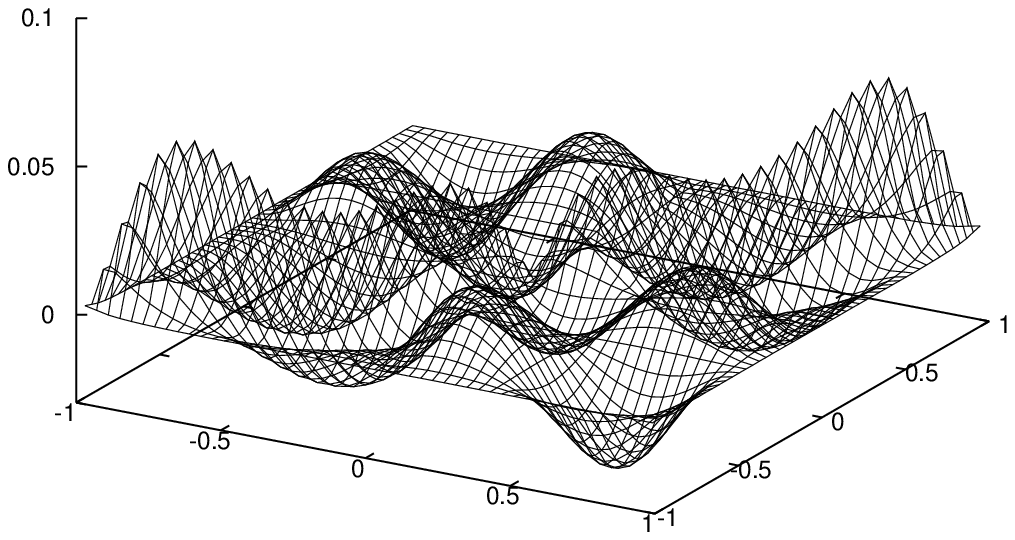}
\caption{MLE-H1 ($\alpha=1.5$, $\kappa=2.5$)}
\label{fig:cov-mle-a1.5-k2.5}
\end{minipage}
\end{center}
\end{figure}

\begin{figure}
\begin{center}
\begin{minipage}{.45\linewidth}
\includegraphics[width=\linewidth]{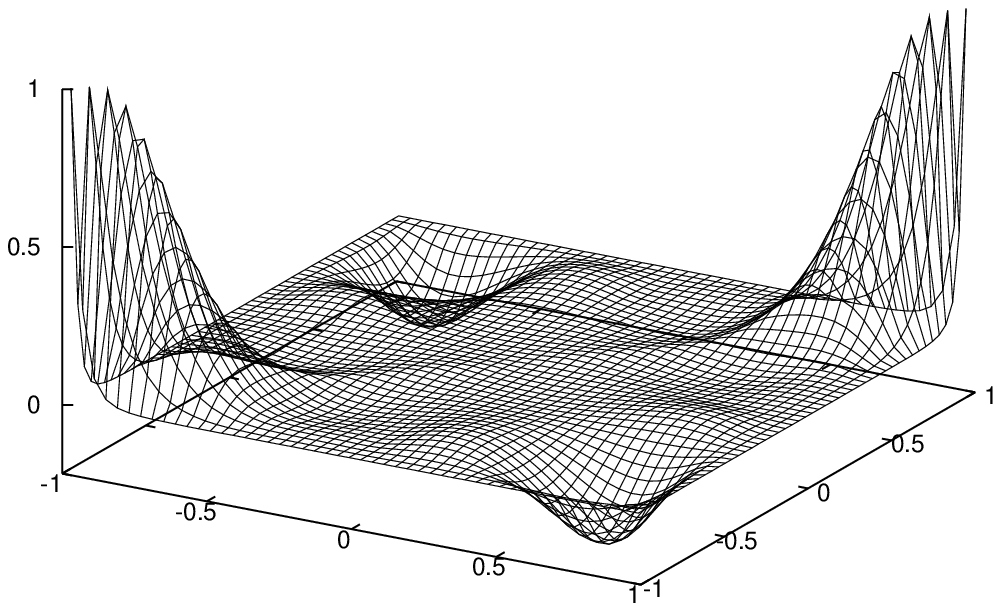}
\caption{MLE-H1 ($\alpha=1.8$, $\kappa=1.0$)}
\label{fig:cov-mle-a1.8-k1.0}
\end{minipage}
\begin{minipage}{.45\linewidth}
\includegraphics[width=\linewidth]{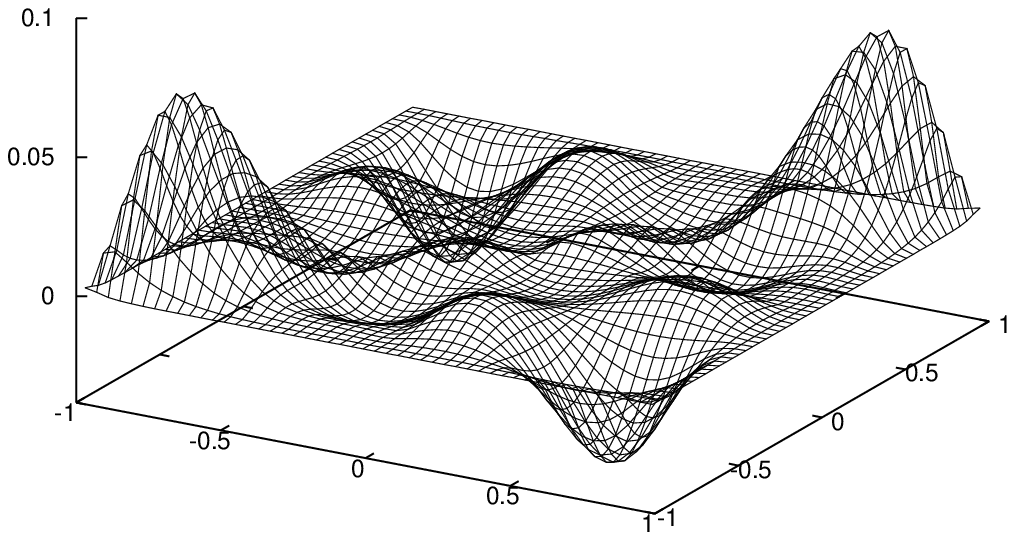}
\caption{MLE-H1 ($\alpha=1.8$, $\kappa=2.5$)}
\label{fig:cov-mle-a1.8-k2.5}
\end{minipage}
\end{center}
\end{figure}

\begin{thm}[Mercer's Theorem, Chapter 5 of Hochstadt (1973)]
\label{mer}
  Let $K(s,t)$ be the kernel of a positive self-adjoint operator on
  $L^2[-1,1]$ and suppose that $K(s,t)$ is continuous in both
  variables. Then
\begin{equation}
  K(s,t)=\sum_{j=1}^{\infty}\frac{1}{\lambda_j}f_j(s)f_j(t),
  \qquad 0 < \lambda_1 \le \lambda_2 \le \cdots \uparrow\infty, 
\label{eq:mercer}
\end{equation}
where $\lambda_j$ is an eigenvalue and $f_j(t)$ is the corresponding
orthonormal eigenfunction 
of the integral equation
\begin{eqnarray}
\lambda \int_{-1}^{1}K(s,t)f(t) dt &=& f(s). \label{eq:int} 
\end{eqnarray} 
The series {\rm (\ref{eq:mercer})} converges uniformly and absolutely to  $K(s,t)$. 
\end{thm}
If $\kappa \le 1$ we need to deal with kernels which are not
continuous at two points $(-1,-1)$ and $(1,1)$. We can see
discontinuity at $(-1,-1)$ and $(1,1)$ in Figures
\ref{fig:cov-mle-a1.0-k1.0}, \ref{fig:cov-mle-a1.5-k1.0} and
\ref{fig:cov-mle-a1.8-k1.0} in the case $\kappa=1$.  On the other hand
there is no discontinuity at these points for $\kappa=2.5$ as shown in
Figures \ref{fig:cov-mle-a1.0-k2.5}, \ref{fig:cov-mle-a1.5-k2.5} and
\ref{fig:cov-mle-a1.8-k2.5}.
However as in Anderson and Darling (1952) the following
version of Mercer's theorem by Hammerstein (1927) is useful.
\begin{thm}
\label{thm:hammerstein}
Suppose that the  covariance function $K(s,t)$ of a Gaussian process  is  
continuous except at $(-1,-1)$ and
$(1,1)$ with $\partial K(s,t)/\partial s$ continuous for 
$|s|,|t| < 1, s\neq t$, and  bounded in $|s| \leq
1-\epsilon$ for every $t\in [-1,1]$
and every $\epsilon >0$.  Then the right hand side of {\rm (\ref{eq:mercer})}
converges uniformly in every domain in the interior of $[-1,1]^2$.
\end{thm}

We apply the above theorems to a continuous
covariance function $K(s,t)$ of a zero mean continuous Gaussian
process $Z(t)$, $-1 < t < 1$, with a finite trace
$\int_{-1}^1 K(t,t)dt < \infty$.
Let $X_1,X_2,\ldots,$ be i.i.d.\  standard normal random
variables. Then the series
$$
Y(t)=\sum_{j=1}^{\infty}\frac{1}{\sqrt{\lambda_j}}f_j(t)X_j
$$
converges in the mean and with probability one for each
$t \in (-1,1)$. Then $Y(t)$ is a Gaussian process with $\mbox{E}Y(t)=0$ and
$\mbox{E}[Y(t)Y(s)]=K(s,t)$. Thus $Y(t)$ defines the same stochastic
process as $Z(t)$. 
Let
\begin{equation}
W^2 = \int_{-1}^{1}Y^{2}(t) dt 
    = \int_{-1}^{1} \biggr\{    
     \sum_{j=1}^{\infty}\frac{1}{\sqrt{\lambda_j}}f_j(t)X_j 
       \biggl\}^{2} dt 
    = \sum_{j=1}^{\infty}\frac{1}{\lambda_j}X_j^{2}. \label{eq:w2y}
\end{equation}
The characteristic function of $W^2$ is given as
$$
\mbox{E}(e^{iuW^2}) 
=  \mbox{E} [  \exp(iu\sum_{j=1}^{\infty }X_{j}^{2}/\lambda_j) ]  
= \prod_{j=1}^{\infty}   \mbox{E}  [  
                      \exp(iuX_{j}^{2}/\lambda_j)  ] \\
= \prod_{j=1}^{\infty}  (1-2iu / \lambda_j)^{- \frac{1}{2}} .
$$
The characteristic function has an alternative expression
$1/\sqrt{D(2it)}$ where $D(\lambda)$ is the associated Fredholm
determinant
$$
 D(\lambda) = \prod_{j=1}^{\infty}   \biggl( 1-\frac{\lambda}{\lambda_j} \biggr).
$$

There are two problems in treating the characteristic function in the
form of the Fredholm determinant. 
One is in the approximation of $D(\lambda)$ itself and the other is in
the Le\'vy's inversion formula. 

In the case of stable distributions the
Fredholm determinant can not be explicitly evaluated and
approximation of $D(\lambda)$ is needed as in the Cauchy case in Matsui and
Takemura (2005).
We approximate $D(\lambda)$ by discretizing the homogeneous integral
equation and approximating eigenvalues of resulting finite system of
linear equations.  
Then the integral equation
(\ref{eq:int}) is approximated by the following finite system of linear
equations
$$
\Tilde{f}= \frac{\lambda}{N}  \Tilde{K}\Tilde{f},
$$
where
$$
\Tilde{K} = \begin{pmatrix}
      K(\xi_1,\xi_1) & \dots & K(\xi_1,\xi_N) \\
       \vdots & & \vdots \\
      K(\xi_N,\xi_1) & \dots & K(\xi_N,\xi_N)
\end{pmatrix}, \qquad 
\Tilde{f} = \left(\begin{array}{c}
f(\xi_{1}) \\
\vdots \\
f(\xi_{N})
\end{array}
\right).
$$
Then the Fredholm determinant is approximated as
\begin{equation*}
\tilde D_N(\lambda) 
= \left|  I - \frac{\lambda}{N}\Tilde{K}   \right| 
= \prod_{j=1}^{N} \left( 1-\frac{\lambda}{\tilde \lambda_j}
\right),
\qquad  0 < \tilde \lambda_1 \le \cdots \le \tilde\lambda_N,  
\label{eq:fd}
\end{equation*}
where $1/\tilde \lambda_j = 1/\tilde \lambda_j(N)$ are the eigenvalues
of $\tilde K/N$.
This method is called a quadrature method and we state 
a version of Theorem 3.4 of Baker~(1977) concerning the convergence of eigenvalues. 

\begin{thm}\label{baker}
Let the eigenvalues $\Tilde{\lambda}_j(N)$ 
be obtained by the quadrature method.
If $K(s,t)$ is positive definite and continuous in $s,t \in [-1,1]$,
$$
\lim_{N \to \infty}\Tilde{\lambda}_j(N)=\lambda_j, 
$$
for each $j$ 
and  
$$
\lim_{N \to \infty} \tilde D_N(\lambda)  = D(\lambda)
$$
for each $\lambda$.
\end{thm}

\begin{rem}
  The covariance functions in {\rm (\ref{eq:transformed-kernal-MLE-H1})} and
  {\rm (\ref{eq:transformed-kernal-MLE-H2})} below do not
  satisfy the conditions of this theorem  if $\kappa \leq 1$.
  However this theorem gives only a sufficient condition for the
  convergence. In our problem the values of $\tilde D_N(\lambda)$ seem to
  converge as we increase $N$ even for the case $\kappa \leq 1$ and the
  resulting value is consistent with our Monte Carlo simulations.
  Therefore in the next section we use the approximation
  of this theorem even for the case $\kappa \leq 1$.  It remains to
  theoretically prove that the approximation is valid for the case
  $\kappa \leq 1$.
\end{rem}

The probability density function of the proposed statistic is given by
inverting the characteristic function $\Phi(t)=1/\sqrt{D(2it)}$.
Since integrand of inversion formula is often wildly oscillating and
converges to 0 slowly, the ordinary numerical integration is difficult
(Section 6.1 of Tanaka (1996)). However we can utilize theory of
complex integration in Slepian (1957) and invert $\Phi(t)$ very
efficiently.  This method of inversion does not seem to be commonly
implemented in statistical computations.

Assuming that the kernel $K(s,t)$ has no multiple eigenvalues and the
number of the eigenvalues are infinite, the density and distribution of $D_\kappa$
are calculated as 
\begin{equation*}
f_D(x)=\frac{1}{\pi}\sum_{k=1}^{\infty}(-1)^k
\int_{\frac{\lambda_{2k-1}}{2}}^{\frac{\lambda_{2k}}{2}}
\frac{e^{-xy}}{\sqrt{\prod_{j=1}^{\infty}\left|1-\frac{2}{\lambda_j} y \right|}} dy,
\end{equation*}
\begin{equation*}
F_D(x)=1-\frac{1}{\pi}\sum_{k=1}^{\infty}(-1)^k
\int_{\frac{\lambda_{2k-1}}{2}}^{\frac{\lambda_{2k}}{2}}
\frac{e^{-xy}}{y\sqrt{\prod_{j=1}^{\infty}\left|1-\frac{2}{\lambda_j} y \right|}} dy.
\end{equation*}
The series are alternating and convenient for checking convergence.
Although the above representations have singularity at each endpoint of integral
range, we can remove the singularity by the following transformation 
$y\mapsto z$ as in Slepian (1957). The
$k$-th integral is transformed by 
\begin{equation}
y=\frac{1}{2}\left(\frac{\lambda_{2k}}{2}-\frac{\lambda_{2k-1}}{2}\right)
\cos \pi z
+\frac{1}{2}\left(\frac{\lambda_{2k}}{2}+\frac{\lambda_{2k-1}}{2}\right), \quad 0
\le z \le 1.\label{eq:transform-y}
\end{equation}
Then
\begin{equation*}
dy=-\pi\sqrt{\left(y-\frac{\lambda_{2k-1}}{2}\right)
\left(\frac{\lambda_{2k}}{2}-y \right)} dz.
\end{equation*}
Hence we obtain the following representations suitable for numerical integration.
\begin{equation}
f_D(x)=\sum_{k=1}^{\infty}(-1)^k \int_0^1
\frac{e^{-xy}\sqrt{\left(\frac{\lambda_{2k}}{2}-y\right)\left(y-\frac{\lambda_{2k-1}}{2}\right)}
}{\sqrt{\prod_{j=1}^{\infty}\left|1-\frac{2}{\lambda_j}y\right|
}} dz,
\end{equation}  
\begin{equation}
F_D(x)=1-\sum_{k=1}^{\infty}(-1)^k \int_0^1
\frac{e^{-xy}\sqrt{\left(\frac{\lambda_{2k}}{2}-y\right)\left(y-\frac{\lambda_{2k-1}}{2}\right)}
}{y\sqrt{\prod_{j=1}^{\infty}\left|1-\frac{2}{\lambda_j}y\right|
}} dz, \label{eq:dist-of-teststatistics}
\end{equation}
where $y$ is given by formula (\ref{eq:transform-y}).

Finally we will make a transformation of
variable and map $\mathbf{R}^2$ into $[-1,1]^2$ in order to satisfy
the finite interval condition of Mercer's theorem. This transformation also is
useful for numerical approximation of eigenvalues.
For deriving the distribution of
$D_\kappa$, we have to incorporate the weight function $e^{-\kappa |t|}$
into the kernel, i.e., we consider the following kernel
$$
\Gamma(s,t)e^{-\frac{\kappa}{2}(|s|+|t|)}.
$$

Now we make the transformation $s\mapsto u$ defined by
$$
 s=-\sgn u\cdot \log(1-|u|), \qquad -1 \le u \le 1.
$$
Then
$$
 ds = \frac{1}{1-|u|}du.
$$
The kernel and the eigenfunctions are transformed as
\begin{eqnarray*}
\Gamma(s,t) &\mapsto&  K(u,v)=
\frac{\Gamma(-\sgn u\cdot \log(1-|u|),-\sgn v\cdot \log(1-|v|))}{\sqrt{(1-|u|)(1-|v|)}}, \\
f_j(s) &\mapsto&  \frac{f_j(-\sgn u\cdot \log(1-|u|))}{\sqrt{1-|u|}}.
\end{eqnarray*}
Eigenvalues of (\ref{eq:int}) do not
change by this transformation and so does Fredholm determinant.
After this transformation, writing  $s,t$ instead of $u,v$ again, 
we have the following kernels on $[-1,1]^2$:

\begin{eqnarray}
\mbox{$H_1$}:  K(s,t) &=& \left\{
e^{-|-\sgn s\cdot \log(1-|s|)+\sgn t\cdot \log(1-|t|)|^\alpha}
-e^{-|\log(1-|s|)|^\alpha-|\log(1-|t|)|^\alpha}\right. \label{eq:transformed-kernal-MLE-H1}\\
&& 
-\Big(
I^{11}\sgn s\cdot \sgn t\cdot \log(1-|s|)\log(1-|t|) \nonumber \\
&& 
+\alpha^2I^{22}|\log(1-|s|)\log(1-|t|)|^\alpha \nonumber \\
&& 
+\alpha I^{23}|\log(1-|s|)\log(1-|t|)|^\alpha
\log|\log(1-|s|)\log(1-|t|)| \nonumber \\
&& 
+I^{33}|\log(1-|s|)\log(1-|t|)|^\alpha
\log|\log(1-|s|)|\cdot\log|\log(1-|t|)| \Big)\nonumber \\
&& 
\left. \quad \times e^{-|\log(1-|s|)|^\alpha-|\log(1-|t|)|^\alpha}\
\right\}\{ (1-|s|)(1-|t|)\}^{\frac{\kappa-1}{2}}.  \nonumber \\
&&          \nonumber  \\
\mbox{$H_2$}:  K(s,t) &=& \left\{
e^{-|-\sgn s\cdot \log(1-|s|)+\sgn t\cdot \log(1-|t|)|^\alpha}
-e^{-|\log(1-|s|)|^\alpha-|\log(1-|t|)|^\alpha}\right. \label{eq:transformed-kernal-MLE-H2}\\
&& 
-\big(
I^{11}\sgn s\cdot \sgn t\cdot
\log(1-|s|)\log(1-|t|)+\alpha^2I^{22}|\log(1-|s|)\log(1-|t|)|^\alpha \big) \nonumber \\
&& \quad 
\left. \times e^{-|\log(1-|s|)|^\alpha-|\log(1-|t|)|^\alpha}\
\right\} \{ (1-|s|)(1-|t|)\}^{\frac{\kappa-1}{2}}.  \nonumber
\end{eqnarray}


\subsection{Numerical approximation of critical values of $D_\kappa$}
We approximate the eigenvalues in (\ref{eq:dist-of-teststatistics}) by
the quadrature method for the kernels (\ref{eq:transformed-kernal-MLE-H1}) and
(\ref{eq:transformed-kernal-MLE-H2}). 800 eigenvalues are calculated by
the above simple algorithm for the case of $\{\kappa=1.0,2.5,5.0,10.0\}$.  
We do not consider $\kappa=0.1$ and $0.5$ since convergence of infinite
integral of $D_{n,\kappa}$ becomes very slow for small weights and we had 
numerical difficulties. 
In Matsui and Takemura (2005) the approximated sum of $\sum_{j=1}^\infty 1/\lambda_j$ did not
converge to $\mbox{E}[D_\kappa]$ first for small $\kappa$. 
We mention that, unlike the Cauchy $\alpha=1$ case,  we did
not observe multiple eigenvalues for other symmetric stable
distributions ($\alpha \neq 1$). 

The infinite series and infinite products in
(\ref{eq:dist-of-teststatistics}) have to be approximated by a finite
sum and finite products. Let $l$ and $m$ $(l<m)$ denote the number of
terms in the sum and the products respectively. Then we can approximate $F_D(x)$
as 
\begin{equation*}
F_D(x) \approx 1-\sum_{k=1}^{l}(-1)^k \int_0^1
\frac{e^{-xy}\sqrt{\left(\frac{\lambda_{2k}}{2}-y\right)\left(y-\frac{\lambda_{2k-1}}{2}\right)}
}{y\sqrt{\prod_{j=1}^{m}\left(1-\frac{2}{\lambda_j}y\right)
}} dz.
\end{equation*}
The series is alternating. Therefore the range of the critical value
can be obtained by substituting lower bound of each positive term and
upper bound of each negative term separately. However deriving
analytical bound of integral of each term of series is difficult.
Hence for accuracy of approximation of $F_D(x)$ we depend on numerical
confirmation. First, we found that finite interval quadrature (QAG) is
very accurate if we set relative error bounds below $10^{-5}$
and the convergence of series is very fast. The first 10 terms of
series are enough to obtain 1\% relative accuracy for most $\kappa$ and $x$
large enough to calculate critical values. Further for $m>300$ the
value of $F_D(x)$ does not change with $m$ for most $\kappa$ and large
$x$.  Finally the approximated percentage points (10\% and 5 \%) 
coincide with simulation results in Section \ref{sec:computational-studies}.

We give Table \ref{tbl:upper-critical-H1} and Table
\ref{tbl:upper-critical-H2} for approximate percentage points of
$D_\kappa$ under hypothesis $H_1$ and $H_2$ respectively. Intervals of
$\alpha \in (0.5,2.0)$ for $H_1$ and $\alpha \in (0.5,2.0]$ are
0.1.  In each table, we set $m=500$ for $\kappa \le 5.0$ and $m=300$
for $\kappa = 10.0$, and set $l=25$ for $\kappa \le 2.5$ and $l=10$
for $\kappa \ge 5.0$. Trial and error indicates that since $D_\kappa$
for large $\kappa$ is very small, we use only accurate large values of
$1/\lambda_j$ among $800$ values.  We also plot the percentage points of each
$D_\kappa$ under $H_1$ in Figures
\ref{fig:upper-teststatistc-1.0}-\ref{fig:upper-teststatistc-10.0}.
The values of percentages are continuous with respect to $\alpha \neq
2$. For large values of $\kappa$, percentage points of small $\alpha$
are large compared to that of large $\alpha$.

\section{Computational studies}
\label{sec:computational-studies}
In this section we give various computational results. Since the exact
finite sample distributions are difficult to obtain, first we
approximate the percentage points of $D_\kappa$ under $H_1$ and $H_2$
respectively by Monte Carlo simulation. Then the power of testing $H_2$ for
the finite sample is evaluated. 

\subsection{Maximum likelihood estimation}
For MLE we maximize likelihood function in parameter space $-\infty
<\mu< \infty$, $\sigma>0$, $0<\alpha\le 2$ by utilizing the first
derivatives of each parameter.  Maximizations are done by the method
based on M.J.D. Powell's TOLMIN from IMSL library.  Although explicit
forms of the density and the derivatives are not available for stable
distributions, the method suggested by Matsui and Takemura (2006)
which improves the original method of Nolan (1977,2001) gives very
accurate approximations.  We use median as the initial value for $\mu$
and for $\sigma$ and $\alpha$ we do grid search and obtain the
parameter values which maximize median inserted log-likelihood among,
say, 2000 points.  Note that although we can set $\alpha=1$ for the
initial value of $\alpha$, the convergence is slow compared with grid
search based initial values.  Base on 1000 Monte Carlo replications,
the values of the estimators and the simulated information matrices
$I(\hat{\theta})=\mbox{Cov}[\hat{\theta}_i,\hat{\theta}_j]^{-1}$ are
given in Table \ref{tbl:MLE-estimation} for
$\alpha \in \{0.8,1.0,1.5,1.8,2.0\}$ and 
the sample size $n \in \{50,100,200\}$. We put true values at the upper row
of each values of $\alpha$ in Table \ref{tbl:MLE-estimation}.  Except
for $\alpha= 2$, simulated values of $I(\hat{\theta})$ coincide
with theoretical values $I(\theta)$ which are given in Matsui and
Takemura (2006).  Though information of $\alpha$ at $\alpha =2.0$ is
$\infty$ and asymptotic normality is not guaranteed (DuMouchel
(1973)), we can estimate $\alpha$ computationally at $\alpha=2$.
Interestingly we observe $\hat{\alpha}=2.0$ for 80\%--90\% of the
cases for $n\ge100$ and we also observe some downward bias.

\subsection{Finite sample critical values of $D_{n,\kappa}$}
We omit the case $\alpha=2$ for both $H_1$ and $H_2$, since there are
many papers concerning testing normality, e.g., Henze and Wagner
(1997), Cs\"{o}rg\H{o} (1986,89) or Naito (1996). Further we
investigate $D_{n,\kappa}$ only $\alpha=1.0,1.5,1.8$ for convenience.
More extensive simulation studies of $D_{n,\kappa}$ for other values of
$\alpha$ are left to our future works.
 
We can compute $D_{n,\kappa}$ by (\ref{eq:teststa}), when the values
of the estimators have converged.  Based on 5000 Monte Carlo
replications, the upper 10 and 5 percentage points of the statistics
$D_{n,\kappa}$, $\kappa \in \{1.0,2.5,5.0,10.0\}$ are tabulated in
Tables \ref{tbl:sim-teststatic-mle-H1-1},
\ref{tbl:sim-teststatic-mle-H1-2}, \ref{tbl:sim-teststatic-mle-H1-3},
\ref{tbl:sim-teststatic-mle-H1-4} for $H_1$ and Table
\ref{tbl:sim-teststatic-mle-H2-1}, \ref{tbl:sim-teststatic-mle-H2-2},
\ref{tbl:sim-teststatic-mle-H2-3}, \ref{tbl:sim-teststatic-mle-H2-4}
for $H_2$. We tabulate simulated values in upper row and asymptotic values
in lower row in box of each value of $\alpha$.  In the tables of $H_1$
when $\alpha$ and $\kappa$ are large, the convergences of
$D_{n,\kappa}$ to $D_{\kappa}$ are slower than other values of
$\alpha$ and $\kappa$. This tendency is also seen in the tables of
$H_2$. This is explained as follows. Since behavior near origin of the
characteristic function reflects behavior of the tail of distribution,
the convergence of the empirical characteristic function near
origin may be faster than that at distant values when the tail of
distribution is heavy.  However we find the values of
$D_{n,\kappa}$ converging $D_{\kappa}$ for all values of $\alpha$ and
$\kappa$ as $n\to \infty$.

\subsection{Analysis of finite sample power}
\label{sec:power}
In this section we examine finite sample power under $H_2$.  For
alternative hypothesis we consider Student's $t$ distribution with $j$
degrees of freedom for $j=1,2,3,4,5,10,\infty$, ($t(1)=C(0,1)$ and
$t(\infty)=N(0,1)$), since stable models are sometimes compared with
Student's $t$ models in empirical applications. We consider only
$\alpha=1.5$ and $\alpha=1.8$ for null distribution of $H_2$ because
the simulation studies are heavy when many values of $\alpha$ are
considered. 
Investigations of other values of $\alpha$ are also left to our future
works.

For the significance levels $\zeta=0.1,0.05$, finite sample power of the
proposed tests are tabulated in Table \ref{tbl:sim-alternative-mle-H2-100-1.5} 
and \ref{tbl:sim-alternative-mle-H2-200-1.8}, based on 1000 Monte
Carlo replications. We summarize our findings. For $\alpha=1.5$
the test with $\kappa=10.0$ has poor power compared with other values of
$\kappa$ and the test with $\kappa=5.0$ is the most powerful. When alternative
hypothesis is $t(3)$ or near $t(3)$ finite sample power is not good. 
For $\alpha=1.8$ the test with $\kappa=10.0$ has also poor power as in
the test for $\alpha=1.5$.  While the test with $\kappa=5.0$ is the
most powerful for heavy tail alternatives, the test with $\kappa=2.5$
is more powerful for the light tail alternatives. These results are
interesting because although the tails of $t(i)$ and $f(x;\alpha)$ are
different, the distributions are not distinguishable well in finite
samples.  

Finally make a remark on how to perform a test of $H_1$. 
The problem is that even the asymptotic null distribution  under $H_1$
depends on the true value of $\alpha$.
In the following remark $D_{n,\kappa}^0$ denotes the 
observed valued of the test statistic for $H_1$ and 
$D_{n,\kappa}(\xi;\alpha)$ denotes the  upper $\xi$ percentage
points of the null distribution of $D_{n,\kappa}$ for $H_1$  when $\alpha$ is the
true characteristic exponent.
\begin{rem}
We can consider several procedures for $H_1$. 
\begin{enumerate}
\item 
$H_1$ is rejected if $D_{n,\kappa}^0 \ge \sup_{\alpha \in (0,2]} D_{n,\kappa}(\xi;\alpha)$,
\item We consider $\alpha \in [a,b]$ where $[a,b] \subset (0,2)$ is a fixed
      range and put $H_1:\alpha \in [a,b]$. Then $H_1$ is rejected if
 $D_{n,\kappa}^0\ge \sup_{\alpha\in[a,b]}\ D_{n,\kappa}(\xi;\alpha) $.
\item We plug in the estimate $\hat{\alpha_n}$ in 
$D_{n,\kappa}(\xi;\hat\alpha_n)$ and $H_1$ is rejected if
$D_{n,\kappa}^0\ge  D_{n,\kappa}(\xi;\hat\alpha_n)$.
\end{enumerate}
Though procedure 1 is logically correct for finite sample, it has the
drawback that the null hypothesis with small true characteristic
exponent $\alpha_0$ may be rejected 
by comparing $D_{n,\kappa}^0$ to percentage points near $\alpha=2$ (See
 Figures 
 \ref{fig:upper-teststatistc-1.0}-\ref{fig:upper-teststatistc-10.0}). 
Therefore procedure 2 and 3 are also worth considering. From the
viewpoint of asymptotic theory we can use procedure 3.  We may use
Procedure 2 considering the standard error in $\hat\alpha_n$.
\end{rem}

\appendix

\section{Proof of Theorem \ref{thm:2.1eise}}
\label{sec:a1}
For simplicity we definite some constants and functions.
$$
w''(x):= d^2 w(x)/{dx}^{2},
\quad 
w_1(t):=w'(|t|)|t|+w(t), 
\quad
w_2(t):=\alpha {|t|}^{\alpha}w(t),
$$
\[
c_1 = \frac{1}{2}\int_{-\infty}^{\infty}e^{-2|t|^{\alpha}}w_1(t)dt, \quad 
c_2 = \int_{-\infty}^{\infty}e^{-2|t|^{\alpha}}w_2(t)dt.
\]
Expanding the estimating equation $Q_\theta(\theta)=(
Q_\mu(\theta),Q_\sigma(\theta),Q_\alpha(\theta))=0$ 
around the true parameter $\theta_0=(0,1,\alpha)$, we have
\[
 Q_\theta(\theta_0)+\frac{\partial Q_\theta(\theta^*)}{\partial \theta}
(\hat{\theta}_n-\theta_0)=0,
\]
where $\theta^\ast_n$ is some value between $\theta_0$ and $\hat{\theta}_n$.
We can write 
\begin{eqnarray}
\sqrt{n}Q_\mu(\theta)
&=& - \frac{1}{\sqrt{n}}\sum_{j=1}^{n}g_1(X_j), \nonumber \\
\sqrt{n}Q_\sigma(\theta)
 &=&  \frac{\sqrt{n}}{2}\left\{ 
   \frac{1}{n^2}
    \sum_{j,k=1}^{n}h_1(X_j,X_k)-\frac{1}{n}
    \sum_{j=1}^{n}2h_2(X_j)\right \},\\
\sqrt{n}Q_\alpha(\theta) 
 &=&
  \frac{1}{\sqrt{n}}\sum_{j=1}^{n}g_2(X_j),  \nonumber
\end{eqnarray}
where 
\begin{eqnarray*}
g_1(x)
  &=&
\int_{-\infty}^{\infty}
     \sin(tx)te^{-|t|^{\alpha}}w(t)dt, \\
g_2(x) 
  &=&
    \int_{-\infty}^{\infty}
     \left\{\cos(tx)-e^{-|t|^{\alpha}}\right\}
     \big({|t|}^{\alpha}\log|t|e^{-|t|^{\alpha}}\big)w(t)dt, \\
h_1(x_1,x_2)
  &=&
\int_{-\infty}^{\infty}
\cos(t(x_1-x_2))w_1(t)dt, \\
h_2(x)
 &=&
 \int_{-\infty}^{\infty}
\cos(tx)e^{-|t|^\alpha}(w_1(t)-w_2(t))dt-c_1+c_2, 
\end{eqnarray*}

$2\sqrt{n}Q_\sigma(\theta)$ can be expressed in the form of a $U$-statistic
$$
2\sqrt{n}Q_\sigma(\theta)=\sqrt{n}\left\{U_n+\frac{h_1(X_1,X_1)}{n}
-\frac{1}{n^2(n-1)}\sum_{j<k}^{n}2h_1(X_j,X_k)\right\}=\sqrt{n}U_n+r_{n},
 \qquad r_{n} \stackrel{P}{\longrightarrow} 0,
$$
where
$$
U_n=
   \frac{2}{n(n-1)}\sum_{1\le j < k \le n}^{n} 
    \left\{h_1(X_j,X_k)-h_2(X_j)-h_2(X_k)\right\}  
=\binom{n}{2}^{-1}
\sum_{1\le j < k \le n}^{n}h(X_j,X_k)  .
$$
By standard argument on $U$-statistics (Chapter 3 of Maesono~(2001),
Chapter 5 of Serfling (1980))
we only need to evaluate 
$$
a(x_1)=\mbox{E}[h(X_1,X_2) \mid X_1=x_1],
$$
since
$$
\sqrt{n}U_n = \frac{1}{\sqrt{n}}\sum_{j=1}^{n}2a(X_j)+r_{n},
 \qquad r_{n} \stackrel{P}{\longrightarrow} 0.
$$
Calculating 
\[
 \mbox{E}[h_1(X_1,X_2)|X_1=x_1]=\int_{-\infty}^{\infty}
\cos(t x_1)e^{-|t|^\alpha}w_1(t)dt
\]
and $\mbox{E}[h_2(X_2)|X_1=x_1]=\mbox{E}[h_2(X_2)]=c_1$,
it can be shown that $a(x_1)$ is written as
$$
a(x_1) 
=\int_{-\infty}^{\infty}\cos(t x_1)e^{-|t|^\alpha}w_2(t)dt-c_2.
$$
After showing the convergence of second derivatives 
$\partial Q_\theta(\theta^*)/\partial \theta
\stackrel{P}{\longrightarrow}A$ we can obtain AL representations,
 \[
\sqrt{n}(\hat{\theta}_n-\theta_0)=-A^{-1} \sqrt{n}Q_\theta(\theta_0).
\]

The proof of $ \partial Q(\theta^*)/\partial \theta
\stackrel{P}{\longrightarrow}A$ is as follows. As before the derivatives
are evaluated at $(\mu,\sigma)=(0,1)$.
Write
{\allowdisplaybreaks
\begin{eqnarray*}
\frac{\partial Q_\mu(\theta)}{\partial \mu}
&=&
\frac{1}{n}\sum_{j=1}^{n} \int_{-\infty}^{\infty}
\cos(tx_j)t^2
e^{-{|t|}^\alpha}w(t) dt,\\
\frac{\partial Q_\sigma(\theta)}{\partial \sigma}
&=&
\frac{1}{2}
\int_{-\infty}^{\infty}
\biggl[
\frac{1}{n^2}\sum_{j,k}\cos(t(x_j-x_k))-\frac{2}{n}\sum_{j=1}^{n}\cos(t x_j)
e^{-{|t|}^\alpha}+e^{-2{|t|}^\alpha}
\biggr]
(w''(  t) |t|^2+2w'(  |t|)|t|)dt \\
&& +
 \frac{1}{n}\sum_{j=1}^{n}\int_{-\infty}^{\infty}
(\cos(t x_j) -e^{-{|t|}^\alpha})
 e^{-{|  t|}^\alpha}|  t|^\alpha 
\left \{(-\alpha^2 |t|^\alpha+ \alpha^2+ \alpha)
   w(t) +2\alpha w'(  t)|t| \right \}dt \\
&&+ \alpha^2 \int_{-\infty}^{\infty}e^{-2{|t|}^\alpha}
{| t|}^{2\alpha} w(t)dt,\\
\frac{\partial Q_\sigma (\theta)}{\partial \alpha}
&=&
\frac{1}{n}\sum_{j=1}^{n}\int_{-\infty}^{\infty}
(\cos(tx_j)-e^{-{|  t|}^\alpha})
e^{-{|  t|}^\alpha}{|  t|}^\alpha 
\left[
\left\{
\log|  t|(-\alpha|  t|^\alpha+\alpha+1)+1\right\}
w(  t) 
 +\log|  t|w'(  t)  |t| 
\right]dt \\
&&+\alpha \int_{-\infty}^{\infty}e^{-2{|t|}^\alpha}
{|t|}^{2\alpha}\log|t| w(t)dt, \\
\frac{\partial Q_\alpha(\theta)}{\partial \alpha}
&=&
\frac{1}{n}\sum_{j=1}^{n}\int_{-\infty}^{\infty}
(\cos(tx_j )-e^{-{|t|}^\alpha})
(1-{|t|}^\alpha)e^{-{|t|}^\alpha}{|t|}^\alpha{(\log|t|)}^2w(t)dt \\
&& +\int_{-\infty}^{\infty}e^{-2{|t|}^\alpha}
{|t|}^{2\alpha}{(\log|t|)}^2w(t)dt.\\
\end{eqnarray*}
}
Making use of $\mbox{E}[ \cos(t X)-e^{-{|t|}^\alpha}]=0$ and
by Fubini's theorem we can calculate their expectations,
\begin{eqnarray*}
\mbox{E}\biggl [\frac{\partial Q_\mu(\theta)}{\partial \mu}  \biggr]
&=& \int_{-\infty}^{\infty}
     e^{-2{|t|}^\alpha}
       t^2 w(t)dt,\\
\mbox{E}\biggl [ \frac{\partial Q_\sigma(\theta)}{\partial \sigma} \biggr]
&=& \frac{1}{2n}
\int_{-\infty}^{\infty}(1-e^{-2{| t|}^\alpha})
(w''( t) |t|^2+2w'( t)|t|)dt \\
&&
+ \alpha^2 \int_{-\infty}^{\infty}e^{-2{|t|}^\alpha}
{|t|}^{2\alpha} w(t)dt,\\
\mbox{E}\biggl [ \frac{\partial Q_\sigma(\theta)}{\partial \alpha} \biggr]
&=&  \alpha \int_{-\infty}^{\infty}
       e^{-2{|t|}^\alpha}
         {|t|}^{2\alpha} \log|t| w(t)dt,\\
\mbox{E}\biggl [\frac{\partial Q_\alpha(\theta)}{\partial \alpha} \biggr]
&=& \int_{-\infty}^{\infty}e^{-2{|t|}^\alpha}
{|t|}^{2\alpha}{(\log|t|)}^2w(t)dt.
\end{eqnarray*}
$\mbox{E}[\partial Q_{\theta_i}(\theta) / \partial \theta_j] =0$ for
other parameters $\theta_i$ since the density is symmetric.
By the weak law of large numbers and continuity of integral about 
parameters we can finish the proof. \qed

\section{Proofs of Section 2.2}
\label{sec:proof-of-section2.2}
\subsection{Proof of Theorem 2.3}
 The idea of proofs is essentially the same as those 
 of G\"urtler and Henze~(2000) based on the original proof of Cs\"org\H{o}~(1983).
 However, since parameter $\alpha$ is additionally estimated and stable distribution
 $f(x;\alpha_0)$ dose not have the $\alpha$-moment
 $\mbox{E}[{|X|}^{\alpha}]$, $\alpha >\alpha_0>0 $,
 we reproduce here the outline of the whole proof.
Note that the kernel $k(s,t)$ is changed from G\"urtler and Henze~(2000).
Before considering Fr\'{e}chet
space $C(\mathbf{R})$, we first assume
the restricted space $C(S)$ of continuous 
functions on a compact subset $S$ with
the supremum norm $\|f\|_{\infty}=\sup_{t\in S}|f(t)|$. 
Letting $k(x,t)=\cos(tx)+i\sin(tx)$, an alternative representation of 
$\hat{Z}_n(t)$ is given by
\begin{equation*}
\begin{split}
\hat{Z}_n(t) 
&= \frac{1}{\sqrt{n}}\sum_{j=1}^{n}\left\{\cos(tx_j) 
  + i\sin(tx_j) - e^{-{|\hat{\sigma}_n t|}^{\hat{\alpha}_n}}
    (\cos(t\hat{\mu_n})+i\sin(t\hat{\mu_n})) \right\}   \\
&= \int k(x,t) d\left\{\sqrt{n}\bigl(F_n(x)
   - F(x;\hat{\theta}_n)\bigl)\right\}. \label{eq:kerneltarns}
\end{split}
\end{equation*}
This is the form of kernel transformed empirical process.
We have to check the condition of $(\mathrm{i})^\ast$ $(\mathrm{ii})^\ast$,
 $(\mathrm{iv})$, $(\mathrm{v})$ and 
$(\mathrm{vi})$ of Cs\"org\H{o}~(1983).
Condition $(\mathrm{i})^\ast$ is satisfied from the definition of
the kernel $k(x,t)=\cos(tx)+i\sin(tx)$.
Condition $(\mathrm{ii})^\ast$ is easy following G\"urtler and Henze~(2000).
For $0 < \epsilon < \alpha_0$,
\begin{eqnarray*}
|\,k(x,s)-k(x,t)\,| &=& \sqrt{{\left(\cos(sx)-\cos(tx)\right)}^2 
                   +{\left(\sin(sx)-\sin(tx)\right)}^2 } \\
                &=& \sqrt{2}\sqrt{1-\cos\left((s-t)x\right)} \\
                &=& 2\,\big|\,\sin((s-t)x/2)\,\big|  \\
                &\le& 4\,|\,s-t\,|^{\epsilon/2}|\,x\,|^{\epsilon/2},                \end{eqnarray*}
and $\mbox{E}[{|X|}^{\epsilon} < \infty]$.
For condition $(\mathrm{iv})$ we can check that the covariance matrix
 $\mbox{E}[l(X) l(X)']$ for MLE and EISE are finite and positive
                   definite from Theorem \ref{thm:2.1mle} and Theorem \ref{thm:2.1eise}. 
Condition $(\mathrm{v})$ for EISE is easy because $l_\theta$ are bounded and
continuously  differentiable from Theorem \ref{thm:2.1eise}. However for MLE we need
some argument since $h_\theta =f_\theta/f$ of Theorem
\ref{thm:2.1mle} has no explicit form. For any compact set $K \in \mathbf{R}$
if $x \in K$, $f(x;\alpha)$ and $f_\theta(x;\alpha)$ are continuously 
differentiable with respect to $x$. Thus $h_\theta$ is exists almost
everywhere and finite $x \in K$. 
To see this we differentiate Fourier inversion formula directly and confirm its integrability.
Note that the density $f(x;\alpha)$ has no singularity. 
However as $x \to \infty$ we have to consider the tail orders of $f(x;\alpha)$ and $f_\theta(x;\alpha)$
going to $0$. We utilize asymptotic expansions of $f(x;\alpha)$ and
$f_\theta(x;\alpha)$ in Matsui and Takemura (2006):
\begin{eqnarray*}
f(x;\alpha)&=& 
  \frac{1}{\pi}\sum^{\infty}_{k=1}\frac{\Gamma(k\alpha+1)}{k!} (-1)^{k-1}
 \sin(\frac{\pi \alpha k }{2}) x^{-k\alpha-1},\\
 f'(x;\alpha) &=& 
\frac{1}{\pi}\sum^{\infty}_{k=1}\frac{\Gamma(\alpha k+2)}{k!} 
 (-1)^k \sin(\frac{\pi \alpha k}{2})x^{-k\alpha-2}, \\
f_\alpha(x;\alpha)&=&
\frac{1}{\pi}\sum^{\infty}_{k=1}\frac{\Gamma'(\alpha k+1)}{(k-1)!}
(-1)^{k-1}\sin\left(\frac{\pi\alpha k}{2}\right)x^{-k\alpha-1}
\\ && 
+\frac{1}{\pi}\sum^{\infty}_{k=1}\frac{\Gamma(\alpha k+1)}{(k-1)!}
(-1)^{k-1}
\left[\frac{\pi}{2}\cos\left(\frac{\pi\alpha k}{2}\right)
   - \log x\; \sin\left(\frac{\pi\alpha k}{2}\right)
     \right]x^{-k\alpha-1}.  
\end{eqnarray*}
{}From expansions and the following relations 
\[
f_\mu(x;\alpha)=-f'(x;\alpha),\quad
f_\sigma(x;\alpha)=-f(x;\alpha)-xf(x;\alpha), 
\]
we can confirm that $l_\theta$'s exist almost everywhere and finite
except for $h_\alpha$ at $x=\infty$ since $h_\alpha(x)=O(\log x)$.
{}From the expansions we see that another condition of $(\mathrm{v})$, i.e.\
\begin{equation}
V_l(u)=\sup_{|x| \le u}\left\{
  \left|l(x;\theta^0)\right|+\left|\frac{\partial}{\partial
      x}l(x;\theta^0)\right|\right\} <\infty,  \qquad \forall u>0,
\end{equation}
is also satisfied. Note that in the proof of Theorem at p.527 of Cs\"org\H{o}
(1983) non-existence of $l_\theta$ at $x=\infty$ is permissible. Therefore
condition $(\mathrm{v})$ holds.

Concerning condition $(\mathrm{vi})$, 
for symmetric stable distributions elements of 
$
\nabla_\theta\Phi(t;\theta)
= (\Phi_\mu(t;\theta), 
    \Phi_\sigma(t;\theta),
\allowbreak    \Phi_\alpha(t;\theta))
$
are written as
\begin{eqnarray*} 
\Phi_\mu(t;\theta)
&=&
\frac{\partial\Phi(t;\alpha,\mu,\sigma)}{\partial \mu}
= \{-\sin(\mu t)+i\cos(\mu t)\}t
e^{-|\sigma t|^{\alpha}}, \\
\Phi_\sigma(t;\theta)
&=&
\frac{\partial\Phi(t;\alpha,\mu,\sigma)}{\partial \sigma}
=-\{\cos(\mu t)+i\sin(\mu t)\}
e^{-|\sigma t|^{\alpha}}{|\sigma t|}^{\alpha}\alpha/{\sigma}, \\
\Phi_\alpha(t;\theta)
&=&
\frac{\partial\Phi(t;\alpha,\mu,\sigma)}{\partial \alpha}
= -\{\cos(\mu t)+i\sin(\mu t)\}
e^{-|\sigma t|^{\alpha}}{|\sigma t|}^{\alpha}\log|\sigma t|.
\end{eqnarray*}
Putting $(\mu,\sigma,\alpha)=\theta_0=(0,1,\alpha)$ we obtain
$$
\nabla_\theta\Phi(t;\theta_0)=\left(\Phi_\mu(t;\alpha),
\Phi_\sigma(t;\alpha), \Phi_\alpha(t;\alpha) \right)=
\big(ite^{-{|t|}^{\alpha}}, -e^{-{|t|}^{\alpha}}{|t|}^{\alpha}\alpha,
-e^{-{|t|}^{\alpha}}{|t|}^{\alpha}\log|t| \big).
$$
Because $\nabla_\theta\Phi(t;\theta)$ is continuous and bounded 
if $ (t,\theta) \in
S \times \Theta_0$, where $\Theta_0$ is a closed parameter space
sufficiently near $\theta_0$, the condition is satisfied.
Therefore the weak convergence of
$\hat{Z}_n(t)$ to a zero mean Gaussian process $Z$ is proved in the
space $(C(S),\|\cdot\|_{\infty})$. Since the compact set $S$ is
arbitrary, 
the space 
$(C(S),\|\cdot\|_{\infty})$ can be extended to Fr\'{e}chet
space $C(\mathbf{R})$ easily.

\medskip
The kernel transform of AL representations are given as follows.
Write  $F_0(x)=F(x,\theta_0)$ for simplicity. \\
1. MLE\ : 
\begin{equation} 
\label{eq:kernel-transform-mle}
\int_{-\infty}^\infty  k(x,s)l_\theta(x)dF_0(x)
=I^{-1}(\theta)\nabla_\theta\Phi(s;\theta).
\end{equation}
2. EISE\ :
\begin{equation}
\label{eq:kernel-transform-eise}
 \int k(x,s)l_\theta(x)dF_0(x)=A^{-1}\int k(x,s)h_\theta(x)dF_0(x),
\end{equation}
where 
\begin{eqnarray*}
\int k(x,s)h_\mu(x)dF_0(x)
&=&
i\int_{-\infty}^{\infty}e^{-{|s-u|}^\alpha-{|u|}^\alpha} 
     u w(u)du,   \\
\int k(x,s)h_\sigma(x)dF_0(x)
&=&
B_\sigma e^{-{|s|}^\alpha}
  -\alpha\int_{-\infty}^{\infty}e^{-{|s-u|}^\alpha-{|u|}^\alpha} 
     {|u|}^\alpha w(u)du,   \\
\int k(x,s)h_\alpha(x)dF_0(x)
&=&
 B_\alpha e^{-{|s|}^\alpha}
  -\int_{-\infty}^{\infty}e^{-{|s-u|}^\alpha-{|u|}^\alpha} 
     {|u|}^\alpha\log |u| w(u)du.
\end{eqnarray*}

Let $\langle\cdot,\cdot\rangle$ denote the standard inner product 
of  $\mathbf{R}^3$ . 
Write 
\begin{eqnarray*}
\Phi(t;\hat{\theta}_n)
&=& \Phi(t;\theta_0) - \big\langle\hat{\theta}_n-\theta_0,\nabla_\theta
\Phi(t,\theta_n^*)\big\rangle  \\
&=& \int k(x,t)dF_0(x) -\big\langle\hat{\theta}_n-\theta_0,
\nabla_\theta
\Phi(t,\theta_n^*)\big\rangle,
\end{eqnarray*}
where $\theta_n^*$ is some value between $\theta_0$ and $\theta_n$.
Note that $\theta_n^*\stackrel{P}{\longrightarrow}\theta_0$. 
Now replace $\sqrt{n}(\hat{\theta}_n-\theta_0)$ by its AL
representations. Then $\hat{Z}_n(t)$ is written as
\begin{align}
\hat{Z}_n(t) 
&= \int
k(x,t)d\left\{\sqrt{n}\bigl(F_n(x)-F_0(x)\bigr)\right\}
  -
  \left\langle\sqrt{n}(\hat{\theta}_n-\theta_0),
   \nabla_\theta\Phi(t;\theta_n^*)\right\rangle 
          \nonumber \\
&= Z_n^*(t)+\Delta_n^{(2)}(t)+\Delta_n^{(3)}(t), \nonumber
\end{align}
where 
\begin{equation*}
Z_n^*(t) := \int
k(x,t)d\left\{\sqrt{n}\bigl(F_n(x)-F_0(x)\big)\right\} 
   - \left\langle\frac{1}{\sqrt{n}}
        \sum_{j=1}^{n}l_\theta(x_j),\nabla_\theta
          \Phi(t;\theta_0)\right\rangle.
\end{equation*}
$Z_n^*$  also converges to $Z$. 
The remainder terms $\Delta_n^{(2)}$ and $\Delta_n^{(3)}$ are
defined by
\begin{eqnarray*}
\Delta_n^{(2)} 
&:=& \left\langle\sqrt{n}(\hat{\theta}_n-\theta_0),
   \nabla_\theta\Phi(t;\theta_0)
  -\nabla_\theta\Phi(t;\theta_n^*)\right\rangle,   \\
\Delta_n^{(3)} 
&:=& -\left\langle\epsilon_n,\nabla_\theta
        \Phi(t;\theta_0)\right\rangle, \qquad \epsilon_n=(r_{n1},r_{n2},r_{n3})'.
\end{eqnarray*} 
These remainder terms satisfy $\displaystyle \sup_{t\in
S}|\Delta_n^{(2)}|\stackrel{P}{\longrightarrow}0$, and $
\displaystyle\sup_{t\in
S}|\Delta_n^{(3)}|\stackrel{P}{\longrightarrow}0 $ by conditions
$(\mathrm{iv})$ and $(\mathrm{vi})$ of Cs\"org\H{o}~(1983). The asymptotic
process $Z$ has an alternative expression
$$
Z(t) = \int k(x,t)dB_{F_0}(x)
    - \left\langle\int l_\theta(x)dB_{F_0}(x),
     \nabla_\theta\Phi(t;\theta_0)\right\rangle , 
$$
where $B_{F_0}(x)$ is the Brownian bridge corresponding to the
distribution function  $F_0$,
having covariance function
$\mbox{E}[ B_{F_0}(s)B_{F_0}(t) ]=F_0(s\wedge t)-F_0(s)F_0(t)$.  
$Z^*$ and $Z$ have the
same covariance function 
\begin{eqnarray}
\label{eq:covariance-general}
\Gamma(s,t) &=& \Phi(s-t;\theta_0)-\Phi(s;\theta_0)\overline{\Phi(t;\theta_0)}+\Phi(s,\theta_0)'
                 \mbox{E}[l_\theta(X_1)l_\theta(X_1)']\overline{\Phi(t,\theta_0)} \\ 
&& \quad  -\left\langle \nabla_\theta\Phi(s;\theta_0),\int \overline{k(x,t)}l_\theta(x)dF_0(x)\right\rangle 
   - \left\langle \overline{\nabla_\theta\Phi(t;\theta_0)},\int k(x,s)l_\theta(x)dF_0(x)\right\rangle .
\nonumber
\end{eqnarray}
Note 
$$
\int k(x,s)\overline{k(x,t)}dF_0(x)
=\int k(x,s-t)dF_0(x)=\Phi(s-t;\theta_0) . 
$$
Evaluating (\ref{eq:covariance-general}) for the case of MLE and EISE
using  (\ref{eq:kernel-transform-mle}) and (\ref{eq:kernel-transform-eise})
proves Theorem \ref{thm:2.3}.
\qed

\subsection{Proof of theorem \ref{thm:efficient}}
We have only to show
\begin{equation}
\label{eq:proof-thm-efficient}
\nabla_\theta\Phi(s;\theta_0)'I^{-1}(\theta)
              \overline{\nabla_\theta\Phi(t;\theta_0)} 
=\left\langle 
\overline{\nabla_\theta\Phi(t;\theta_0)},\int k(x,s)l_E(x)dF_0(x)\right\rangle .
\end{equation}
Because AL representations can be written by 
$I^{-1}(\theta) \partial \log f(x;\theta_0)/\partial \theta$, their kernel
transformations are 
\begin{eqnarray*}
\int k(x,s)l_E(x)dF_0(x) &=& I^{-1}(\theta)\int
 k(x,s)\frac{\partial \log f(x;\theta_0)}{\partial \theta} dF_0(x). \\
 &=& 
I^{-1} (\theta)\int k(x,s)\frac{1}{f(x;\theta_0)} \frac{\partial f(x;\theta_0)}{\partial \theta}\ f(x;\theta_0)dx \\
&=&
I^{-1} (\theta) \frac{1}{2\pi} \int e^{isx}\frac{\partial f(x;\theta_0)}{\partial \theta}  dx \\
&=& I^{-1} (\theta) \nabla_\theta\Phi(s;\theta_0).
\end{eqnarray*}
Since both sides of the formula (\ref{eq:proof-thm-efficient}) are scalars
the proof is over.  
\qed

\subsection{Proof of corollary 2.2.1}
Let $\alpha=1$.  Then $I^{-1}({\theta})$ 
is explicitly written as
\[
I^{-1}({\theta})=
\left(
\begin{array}{ccc}
2  & 0                                 &  0                           \\
0 & 2+\frac{12}{\pi^2}{(\gamma+\log2-1)}^2 & \frac{12}{\pi^2}(\gamma+\log2-1)  \\
0  & \frac{12}{\pi^2}(\gamma+\log2-1) & \frac{12}{\pi^2}
\end{array}
\right).
\]
Letting $\alpha=1$ in (\ref{eq:mle-covariance}) and replacing
$I^{-1}({\theta})$ by the above matrix we can prove the corollary.
\qed

\begin{table}[h]
\caption[percentage]{Upper $\xi$ percentage points of
 $D_{\kappa}$ under $H_1$}
\begin{center}
\begin{tabular}{|c||c|c|c|c|c|}    \hline
$\alpha$&$\xi$$\setminus\kappa$  
& 1.0 & 2.5 & 5.0 & 10.0  \\ \hline
1.9 & 0.1  & 1.097 & 0.1030 & 0.00759 & 0.000973 \\ \cline{2-6}
 & 0.05    & 1.352 & 0.1338 & 0.01000 & 0.001267 \\ \hline
1.8 & 0.1  & 1.037 & 0.0963 & 0.00977 & 0.002283 \\ \cline{2-6}
 & 0.05    & 1.273 & 0.1240 & 0.01257 & 0.002984 \\ \hline
1.7 & 0.1  & 0.991 & 0.0958 & 0.01375 & 0.003974 \\ \cline{2-6}
 & 0.05    & 1.211 & 0.1217 & 0.01754 & 0.005183 \\ \hline
1.6 & 0.1  & 0.957 & 0.1007 & 0.01921 & 0.006077 \\ \cline{2-6}
 & 0.05    & 1.161 & 0.1258 & 0.02431 & 0.007895 \\ \hline
1.5 & 0.1  & 0.933 & 0.1108 & 0.02615 & 0.008650 \\ \cline{2-6}
 & 0.05    & 1.122 & 0.1361 & 0.03280 & 0.011180 \\ \hline
1.4 & 0.1  & 0.918 & 0.1260 & 0.03469 & 0.011770 \\ \cline{2-6}
 & 0.05    & 1.094 & 0.1527 & 0.04313 & 0.015120 \\ \hline
1.3 & 0.1  & 0.915 & 0.1462 & 0.04503 & 0.015530 \\ \cline{2-6}
 & 0.05    & 1.078 & 0.1754 & 0.05553 & 0.019813 \\ \hline
1.2 & 0.1  & 0.925 & 0.1717 & 0.05740 & 0.020042 \\ \cline{2-6}
 & 0.05    & 1.078 & 0.2041 & 0.07025 & 0.025370 \\ \hline
1.1 & 0.1  & 0.948 & 0.2027 & 0.07206 & 0.025434 \\ \cline{2-6}
 & 0.05    & 1.095 & 0.2390 & 0.08755 & 0.031916 \\ \hline
1.0 & 0.1  & 0.988 & 0.2395 & 0.08923 & 0.031846 \\ \cline{2-6}
 & 0.05    & 1.130 & 0.2804 & 0.10765 & 0.039574 \\ \hline
0.9 & 0.1  & 1.044 & 0.2824 & 0.10903 & 0.039422 \\ \cline{2-6}
 & 0.05    & 1.186 & 0.3288 & 0.13063 & 0.048454 \\ \hline
0.8 & 0.1  & 1.118 & 0.3315 & 0.13145 & 0.048287 \\ \cline{2-6}
 & 0.05    & 1.262 & 0.3840 & 0.15627 & 0.058625 \\ \hline
0.7 & 0.1  & 1.213 & 0.3855 & 0.15611 & 0.058531 \\ \cline{2-6}
 & 0.05    & 1.362 & 0.4446 & 0.18397 & 0.070081 \\ \hline
0.6 & 0.1  & 1.325 & 0.4413 & 0.18217 & 0.070175 \\ \cline{2-6}
 & 0.05    & 1.482 & 0.5065 & 0.21239 & 0.082722 \\ \hline
0.5 & 0.1  & 1.441 & 0.4928 & 0.20834 & 0.083189 \\ \cline{2-6}
 & 0.05    & 1.609 & 0.5615 & 0.23966 & 0.096393 \\ \hline
\end{tabular}
\label{tbl:upper-critical-H1}
\end{center}
\end{table}

\begin{table}
\caption[percentage]{Upper $\xi$ percentage points of
 $D_{\kappa}$ under $H_2$}
\begin{center}
\begin{tabular}{|c||c|c|c|c|c|}    \hline
$\alpha$&$\xi$$\setminus\kappa$  
& 1.0 & 2.5 & 5.0 & 10.0  \\ \hline
2.0 & 0.1  & 1.216 & 0.1258 & 0.00881 & 0.000241 \\ \cline{2-6}
 & 0.05    & 1.499 & 0.1622 & 0.01177 & 0.000335 \\ \hline
1.9 & 0.1  & 1.150 & 0.1129 & 0.00921 & 0.00142 \\ \cline{2-6}
 & 0.05    & 1.413 & 0.1444 & 0.01164 & 0.00179 \\ \hline
1.8 & 0.1  & 1.110 & 0.1111 & 0.01354 & 0.00329 \\ \cline{2-6}
 & 0.05    & 1.357 & 0.1398 & 0.01679 & 0.00416 \\ \hline
1.7 & 0.1  & 1.080 & 0.1157 & 0.01984 & 0.00566 \\ \cline{2-6}
 & 0.05    & 1.313 & 0.1431 & 0.02457 & 0.00713 \\ \hline
1.6 & 0.1  & 1.058 & 0.1256 & 0.02773 & 0.00858 \\ \cline{2-6}
 & 0.05    & 1.277 & 0.1532 & 0.03426 & 0.01076 \\ \hline
1.5 & 0.1  & 1.044 & 0.1404 & 0.03721 & 0.01211 \\ \cline{2-6}
 & 0.05    & 1.249 & 0.1697 & 0.04578 & 0.01514 \\ \hline
1.4 & 0.1  & 1.037 & 0.1599 & 0.04840 & 0.01636 \\ \cline{2-6}
 & 0.05    & 1.231 & 0.1919 & 0.05927 & 0.02037 \\ \hline
1.3 & 0.1  & 1.039 & 0.1840 & 0.06148 & 0.02144 \\ \cline{2-6}
 & 0.05    & 1.222 & 0.2195 & 0.07493 & 0.02660 \\ \hline
1.2 & 0.1  & 1.051 & 0.2127 & 0.07670 & 0.02748 \\ \cline{2-6}
 & 0.05    & 1.225 & 0.2524 & 0.09300 & 0.03396 \\ \hline
1.1 & 0.1  & 1.074 & 0.2465 & 0.09428 & 0.03464 \\ \cline{2-6}
 & 0.05    & 1.242 & 0.2907 & 0.11376 & 0.04262 \\ \hline
1.0 & 0.1  & 1.111 & 0.2862 & 0.11445 & 0.04307 \\ \cline{2-6}
 & 0.05    & 1.276 & 0.3356 & 0.13742 & 0.05273 \\ \hline
0.9 & 0.1  & 1.159 & 0.3305 & 0.13723 & 0.05291 \\ \cline{2-6}
 & 0.05    & 1.322 & 0.3855 & 0.16395 & 0.06443 \\ \hline
0.8 & 0.1  & 1.226 & 0.3811 & 0.16263 & 0.06426 \\ \cline{2-6}
 & 0.05    & 1.389 & 0.4424 & 0.19326 & 0.07779 \\ \hline
0.7 & 0.1  & 1.310 & 0.4365 & 0.19018 & 0.07714 \\ \cline{2-6}
 & 0.05    & 1.476 & 0.5045 & 0.22466 & 0.09278 \\ \hline
0.6 & 0.1  & 1.412 & 0.4936 & 0.21890 & 0.09144 \\ \cline{2-6}
 & 0.05    & 1.583 & 0.5678 & 0.25681 & 0.10919 \\ \hline
0.5 & 0.1  & 1.517 & 0.5464 & 0.24726 & 0.10689 \\ \cline{2-6}
 & 0.05    & 1.696 & 0.6250 & 0.28773 & 0.12661 \\ \hline
\end{tabular}
\label{tbl:upper-critical-H2}
\end{center}
\end{table}

\begin{table}
\caption[percentage]{Simulation results of symmetric stable distributions
 (1000 iterations)}
\begin{center}
\begin{tabular}{c|ccccccccc}    \hline
n & $\hat{\alpha} $ & $\hat{\mu}$ & $\hat{\sigma}$ & $\hat{I}_{11}$ & $\hat{I}_{22}$ & $\hat{I}_{33}$ 
 &$\hat{I}_{12}$ &$\hat{I}_{13}$ & $\hat{I}_{23}$ \\ \hline
&2.0   &0& 1 &0.5 & 2.0 & $ \infty $     & 0 &0 & $\ast $ \\
50 &1.976 & 0.00014 & 0.977  & 0.607  & 1.875 & 5.647 & -0.029 & -0.010 & -0.828 \\ 
100&1.990 & 0.00017 & 0.975  & 0.868  & 1.402 & 9.445 &  0.032 & -0.057 & -0.594\\ 
200&1.994 & 0.00094 & 0.977  & 0.784  & 1.239 & 12.77 & -0.009 &  0.078&-0.539 \\ \hline
&1.8  &0&1  & 0.4552 & 1.3898 & 0.5937 &0 &0 & -0.3138  \\
50  &1.818 & 0.0012 &  0.991 & 0.487  & 1.231 & 0.676  & -0.033 &-0.026 & -0.267 \\ 
100 &1.822 & 0.0033 & 1.002 & 0.482  & 1.356 & 0.584  & 0.005  &-0.016 &-0.340 \\ 
200 &1.810 &-0.0000& 1.000  & 0.450 & 1.399 &  0.603  & -0.007& -0.003&-0.323 \\ \hline
&1.5   &   0       &1  & 0.4281 & 0.9556 & 0.4737 &0&0 & -0.2174  \\
50&1.548 & -0.0022 & 1.012 & 0.3161 & 0.5796 & 0.4252 & 0.005 &-0.009 & -0.1927 \\ 
100&1.524 & -0.0000 & 1.000 & 0.3914 & 0.9138 & 0.4278 & 0.001 &-0.010 & -0.2291 \\ 
200&1.510 & -0.0003 & 1.000 & 0.4028 & 0.9474 & 0.4683&-0.018&-0.006&-0.2229 \\ \hline
&1.0 &0& 1 & 0.5    & 0.5    & 0.8590 &0&0 & -0.1352  \\
50 &1.026 & -0.0029 & 0.996 & 0.4243 & 0.4877 & 0.670  &-0.004 & 0.0248&-0.1779 \\ 
100 &1.001 & -0.0041 & 0.988 & 0.4438 & 0.4205 & 0.746  &0.007 &0.003&-0.1527 \\ 
200 &1.006& -0.0039 &1.001 & 0.4929 & 0.5013 & 0.845 &-0.003&-0.0251&-0.1534 \\ \hline
&0.8  &0&1 & 0.6800 & 0.3586 & 1.3928 &0&0 & -0.0913  \\
50  &0.815 & -0.0016 & 1.005 & 0.5434 &0.3243 &1.111 &-0.01432 &0.0015 &-0.083 \\ 
100 &0.811 & -0.0001 & 1.003 & 0.6015 &0.3459 &1.171 &0.00435&-0.0313&-0.097 \\ 
200 &0.805 & -0.0009 & 1.000 & 0.6232 &0.3708 &1.303&0.01785&-0.0276&-0.086  \\ \hline
\end{tabular}
\label{tbl:MLE-estimation}
\end{center}
\end{table}

\begin{figure}[h]
\begin{center}
\begin{minipage}{.45\linewidth}
\includegraphics[width=\linewidth]{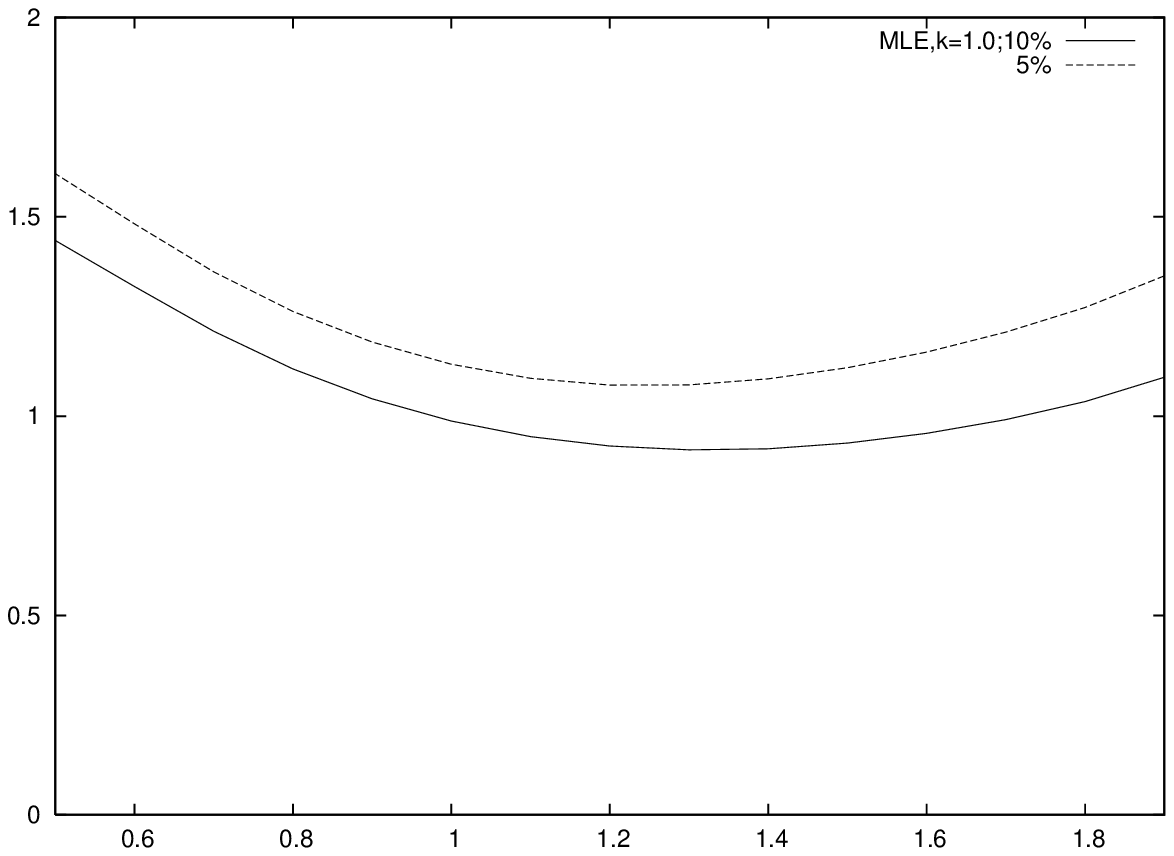}
\caption{Upper quantiles ($\kappa=1.0$, $H_1$)}
\label{fig:upper-teststatistc-1.0}
\end{minipage}
\begin{minipage}{.45\linewidth}
\includegraphics[width=\linewidth]{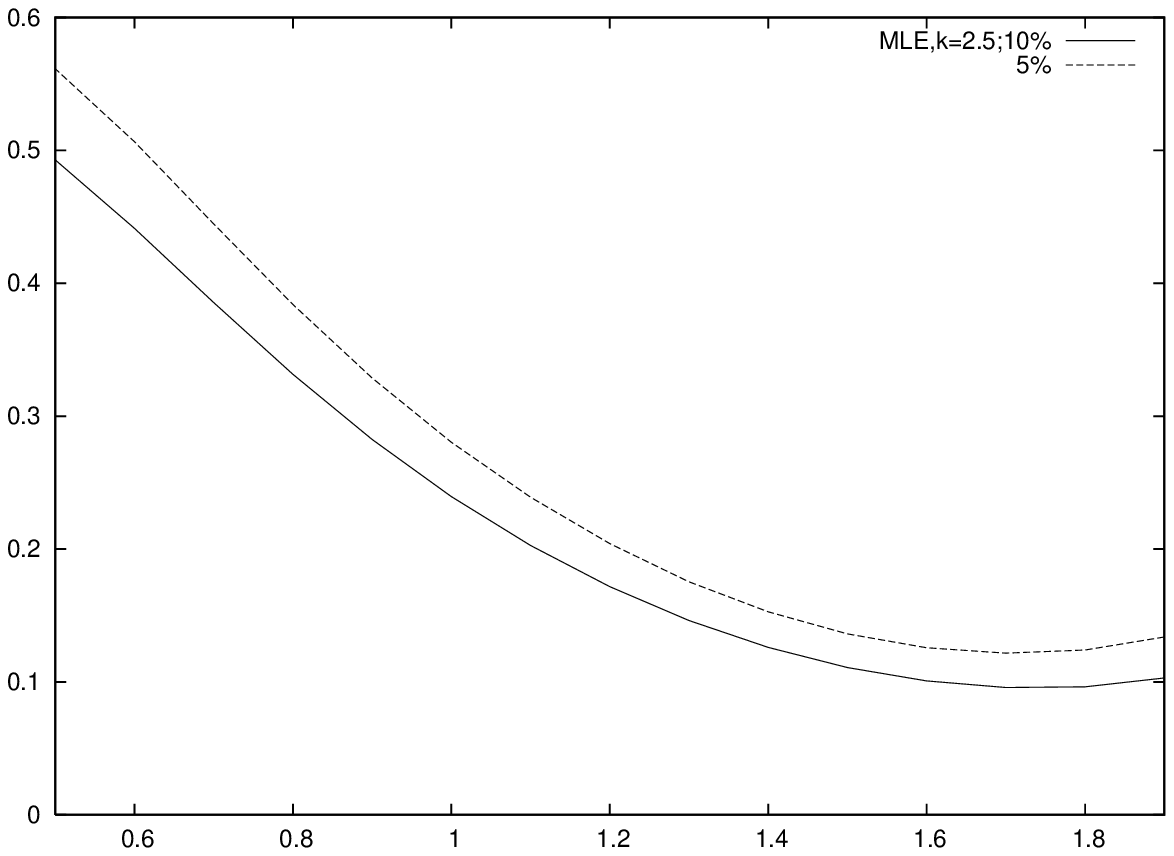}
\caption{Upper quantiles ($\kappa=2.5$, $H_1$)}\label{fig:upper-teststatistc-2.5}
\end{minipage}
\end{center}
\end{figure}

\begin{figure}[h]
\begin{center}
\begin{minipage}{.45\linewidth}
\includegraphics[width=\linewidth]{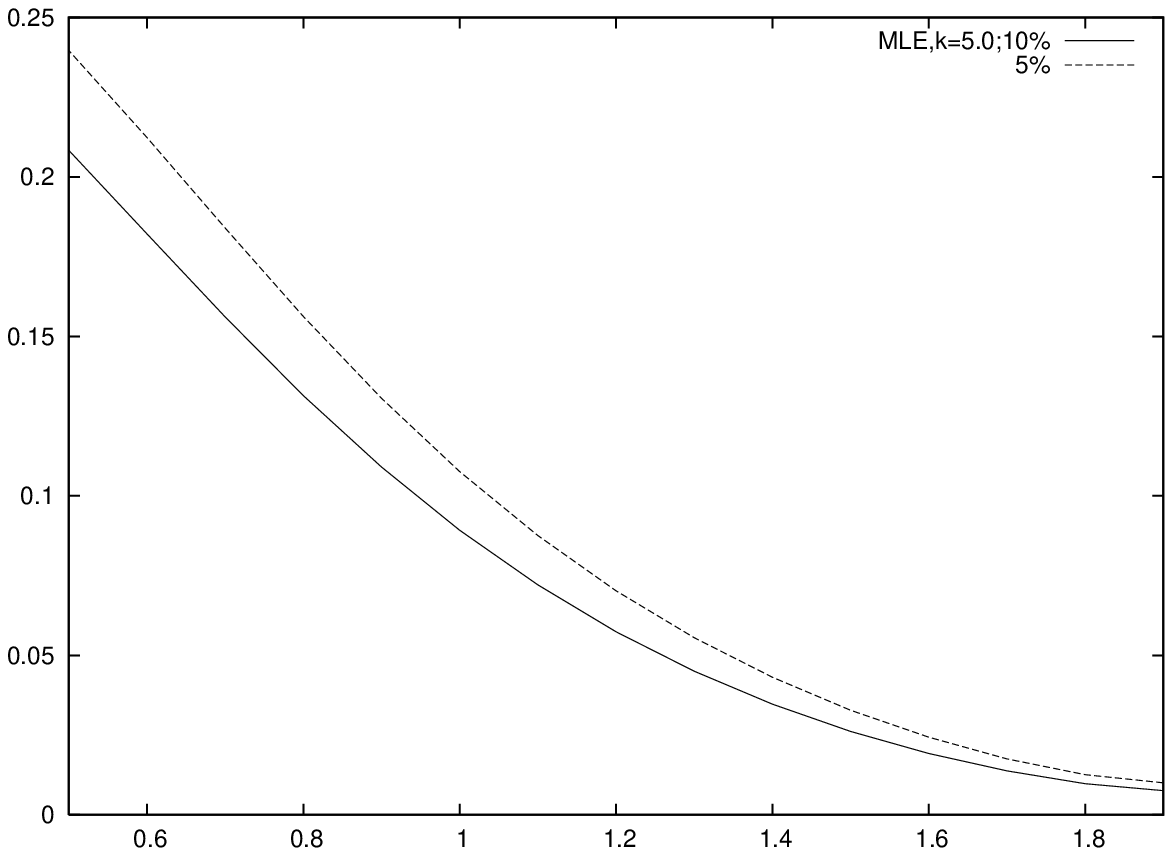}
\caption{Upper quantiles ($\kappa=5.0$, $H_1$)}\label{fig:upper-teststatistc-5.0}
\end{minipage}
\begin{minipage}{.45\linewidth}
\includegraphics[width=\linewidth]{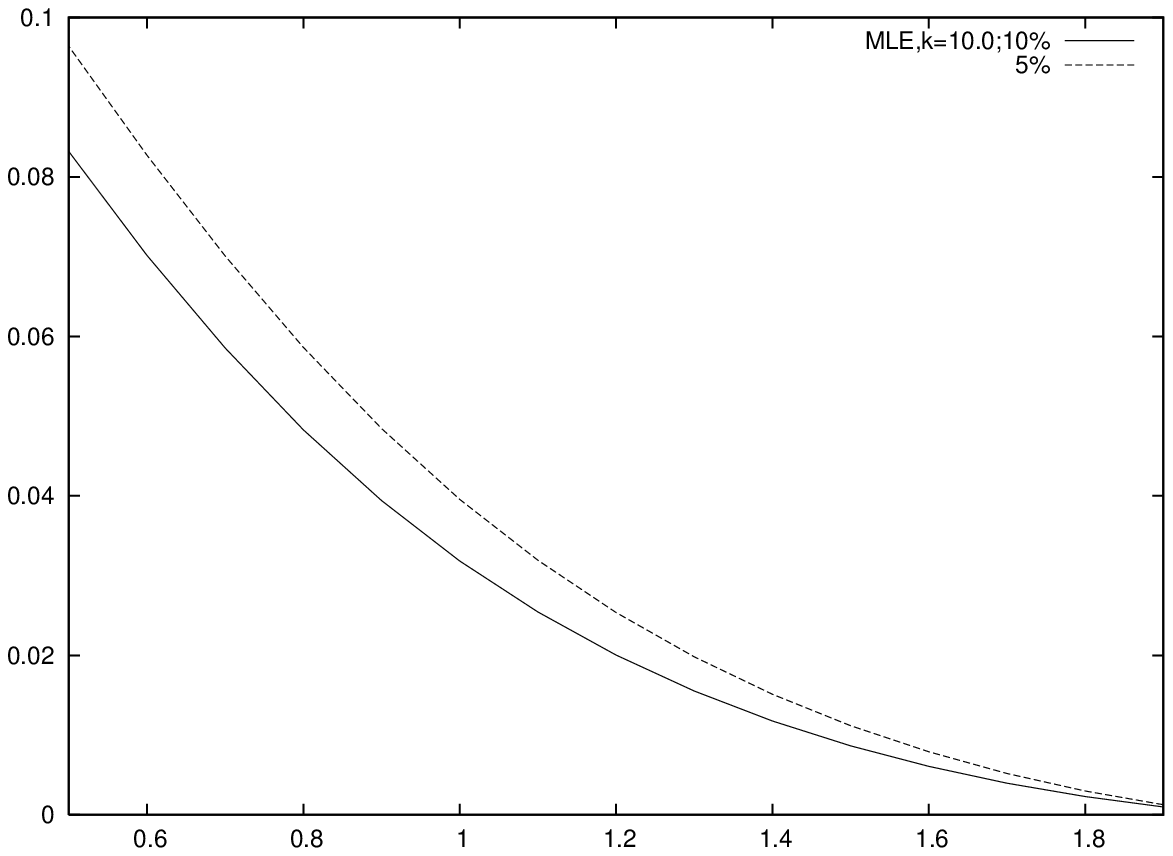}
\caption{Upper quantiles ($\kappa=10.0$, $H_1$)}\label{fig:upper-teststatistc-10.0}
\end{minipage}
\end{center}
\end{figure}

\begin{table}[h]
\caption[percentage]{Upper 10 percentage points of $D_{100,\kappa}$
 under $H_1$}
\begin{center}
\begin{tabular}{ccccc}    \hline
$ \alpha \setminus\kappa$   & 1.0 & 2.5 & 5.0 & 10.0  \\ \hline
1.8  & 1.037 & 0.1037 & 0.01271 & 0.00326 \\
     & 1.037 & 0.0963 & 0.00977 & 0.00228 \\ \hline 
1.5  & 0.953 & 0.1213 & 0.03061 & 0.00996 \\
     & 0.933 & 0.1108 & 0.02615 & 0.00865 \\ \hline 
1.0  & 1.032 & 0.2541 & 0.09450 & 0.03336 \\
     & 0.988 & 0.2395 & 0.08923 & 0.03185
\end{tabular}\label{tbl:sim-teststatic-mle-H1-1}
\end{center}
\end{table}
\begin{table}
\caption[percentage]{Upper 5 percentage points of $D_{100,\kappa}$
 under $H_1$}
\begin{center}
\begin{tabular}{ccccc}    \hline
$ \alpha \setminus\kappa$   & 1.0 & 2.5 & 5.0 & 10.0  \\ \hline
1.8  & 1.279 & 0.1334 & 0.01681 & 0.00451 \\
     & 1.273 & 0.1240 & 0.01257 & 0.00298 \\ \hline
1.5  & 1.141 & 0.1504 & 0.03968 & 0.01348 \\
     & 1.122 & 0.1361 & 0.03280 & 0.01118 \\ \hline
1.0  & 1.216 & 0.3059 & 0.11400 & 0.04230 \\
     & 1.130 & 0.2804 & 0.10765 & 0.03957
\end{tabular}\label{tbl:sim-teststatic-mle-H1-2}
\end{center}
\end{table}

\begin{table}
\caption[percentage]{Upper 10 percentage points of $D_{200,\kappa}$
 under $H_1$}
\begin{center}
\begin{tabular}{ccccc}    \hline
$ \alpha \setminus\kappa$   & 1.0 & 2.5 & 5.0 & 10.0  \\ \hline
1.8  & 1.029 & 0.1003 & 0.01119 & 0.00279 \\
     & 1.037 & 0.0963 & 0.00977 & 0.00228 \\ \hline
1.5  & 0.929 & 0.1145 & 0.02807 & 0.00926 \\
     & 0.933 & 0.1108 & 0.02615 & 0.00865 \\ \hline
1.0  & 1.006 & 0.2462 & 0.09072 & 0.03233 \\
     & 0.988 & 0.2395 & 0.08923 & 0.03185
\end{tabular}\label{tbl:sim-teststatic-mle-H1-3}
\end{center}
\end{table}
\begin{table}
\caption[percentage]{Upper 5 percentage points of $D_{200,\kappa}$ under
 $H_1$}
\begin{center}
\begin{tabular}{ccccc}    \hline
$ \alpha \setminus\kappa$   & 1.0 & 2.5 & 5.0 & 10.0  \\ \hline
1.8  & 1.290 & 0.1292 & 0.01467 & 0.00386 \\
     & 1.273 & 0.1240 & 0.01257 & 0.00298 \\ \hline 
1.5  & 1.125 & 0.1419 & 0.03575 & 0.01222 \\
     & 1.122 & 0.1361 & 0.03280 & 0.01118 \\ \hline
1.0  & 1.161 & 0.2887 & 0.10938 & 0.04015 \\
     & 1.130 & 0.2804 & 0.10765 & 0.03957
\end{tabular}\label{tbl:sim-teststatic-mle-H1-4}
\end{center}
\end{table}

\begin{table}
\caption[percentage]{Upper 10 percentage points of $D_{100,\kappa}$
 under $H_2$}
\begin{center}
\begin{tabular}{ccccc}    \hline
$ \alpha \setminus\kappa$   & 1.0 & 2.5 & 5.0 & 10.0  \\ \hline
1.8  &1.100 & 0.1126 & 0.01415 & 0.00356 \\
     &1.110 & 0.1111 & 0.01354 & 0.00329 \\ \hline
1.5  &1.030 & 0.1403 & 0.03709 & 0.01202 \\
     &1.044 & 0.1404 & 0.03721 & 0.01211
\end{tabular}\label{tbl:sim-teststatic-mle-H2-1}
\end{center}
\end{table}
\begin{table}
\caption[percentage]{Upper 5 percentage points of $D_{100,\kappa}$
 under $H_2$}
\begin{center}
\begin{tabular}{ccccc}    \hline
$ \alpha \setminus\kappa$   & 1.0 & 2.5 & 5.0 & 10.0  \\ \hline
1.8  &1.333  & 0.1397 & 0.01878 & 0.00533 \\
     &1.357  & 0.1398 & 0.01679 & 0.00416 \\ \hline
1.5  &1.220  & 0.1690 & 0.04710 & 0.01590 \\
     &1.249  & 0.1697 & 0.04578 & 0.01514 
\end{tabular}\label{tbl:sim-teststatic-mle-H2-2}
\end{center}
\end{table}

\begin{table}
\caption[percentage]{Upper 10 percentage points of $D_{200,\kappa}$
 under $H_2$}
\begin{center}
\begin{tabular}{ccccc}    \hline
$ \alpha \setminus\kappa$   & 1.0 & 2.5 & 5.0 & 10.0  \\ \hline
1.8  &1.077 &0.1108 & 0.01369  & 0.00332 \\
     &1.110 &0.1111 & 0.01354  & 0.00329 \\ \hline
1.5  &1.037 &0.1414 & 0.03754 & 0.01204 \\
     &1.044 &0.1404 & 0.03721 & 0.01211 
\end{tabular}\label{tbl:sim-teststatic-mle-H2-3}
\end{center}
\end{table}
\begin{table}
\caption[percentage]{Upper 5 percentage~points of $D_{200,\kappa}$ under
 $H_2$}
\begin{center}
\begin{tabular}{ccccc}    \hline
$ \alpha \setminus\kappa$   & 1.0 & 2.5 & 5.0 & 10.0  \\ \hline
1.8  &1.368 & 0.1409 & 0.01818 & 0.00470 \\
     &1.357 & 0.1398 & 0.01679 & 0.00416 \\ \hline
1.5  &1.244 & 0.1695 & 0.04620 & 0.01561 \\
     &1.249 & 0.1697 & 0.04578 & 0.01514 
\end{tabular}\label{tbl:sim-teststatic-mle-H2-4}
\end{center}
\end{table}

\begin{table}
\begin{center}
\caption[percentage]{Power of $D_{n,\kappa}$ under $H_2$ ($\alpha=1.5$,\ significance
  levels~$\xi = 0.1,\ 0.05,\ n=100$)}
\ \\
\begin{tabular}{|l||cccc||cccc|}   
\hline
\multicolumn{1}{|c||}{$\xi$}  &
\multicolumn{4}{c||}{$0.1$} &
\multicolumn{4}{c|}{$0.05$} \\
 \hline
\multicolumn{1}{|c||}{$\kappa$}   & 1.0 & 2.5 & 5.0 & 10.0 & 1.0 & 2.5 & 5.0 & 10.0 \\ \hline
$N(0,2)$   & 44 & 68 & 76 & 45 & 30 & 49 & 45 & 2  \\
$t(1)  $   & 82 & 93 & 96 & 95 & 75 & 90 & 92 & 92 \\
$t(2)  $   & 16 & 20 & 17 & 15 &  9 & 12 & 10 & 8  \\
$t(3)  $   & 11 & 10 &  5 &  3 &  6 &  6 &  2 & 1  \\
$t(4)  $   & 12 & 12 &  9 &  3 &  7 &  6 &  3 & 0  \\
$t(5)  $   & 14 & 15 & 15 &  5 & 8  & 9  & 5  & 0  \\
$t(10) $   & 26 & 40 & 42 & 15 & 14 & 22 & 16 & 0  \\
\hline
\end{tabular}\label{tbl:sim-alternative-mle-H2-100-1.5}
\end{center}
\end{table} 

\begin{table}[h]
\begin{center}
\caption[percentage]{Power of $D_{n,\kappa}$ under $H_2$ ($\alpha=1.5$,\ significance
  levels~$\xi = 0.1,\ 0.05,\ n=200$)}
\ \\
\begin{tabular}{|l||cccc||cccc|}   
\hline
\multicolumn{1}{|c||}{$\xi$}  &
\multicolumn{4}{c||}{$0.1$} &
\multicolumn{4}{c|}{$0.05$} \\
 \hline
\multicolumn{1}{|c||}{$\kappa$}   & 1.0 & 2.5 & 5.0 & 10.0 & 1.0 & 2.5 & 5.0 & 10.0 \\ \hline
$N(0,2)$   & 79 & 99 &100 & 100& 63 & 96 & 99 & 98  \\
$t(1)  $   & 98 & 99 &100 & 100& 96 & 99 &100 &100  \\
$t(2)  $   & 24 & 26 & 21 & 16 & 16 & 18 & 14 & 10  \\
$t(3)  $   & 12 & 10 & 10 &  9 &  6 &  6 &  5 &  3  \\
$t(4)  $   & 16 & 24 & 32 & 27 &  8 & 12 & 16 &  8  \\
$t(5)  $   & 23 & 40 & 50 & 47 & 12 & 26 & 33 & 19  \\
$t(10) $   & 48 & 83 & 93 & 92 & 34 & 66 & 84 & 68  \\
\hline
\end{tabular}\label{tbl:sim-alternative-mle-H2-200-1.5}
\end{center}
\end{table} 

\begin{table}
\begin{center}
\caption[percentage]{Power of $D_{n,\kappa}$ under $H_2$ ($\alpha=1.8$,\ significance
  levels~$\xi = 0.1,\ 0.05,\ n=100$)}
\ \\
\begin{tabular}{|l||cccc||cccc|}   
\hline
\multicolumn{1}{|c||}{$\xi$}  &
\multicolumn{4}{c||}{$0.1$} &
\multicolumn{4}{c|}{$0.05$} \\
 \hline
\multicolumn{1}{|c||}{$\kappa$}   & 1.0 & 2.5 & 5.0 & 10.0 & 1.0 & 2.5 & 5.0 & 10.0 \\ \hline
$N(0,2)$   & 14 & 14 & 3 & 0 & 8 & 8 & 1 & 0  \\
$t(1)  $   & 98 & 100 & 100 & 100 & 96 & 100 & 100 & 100 \\
$t(2)  $   & 55 & 72 & 77 & 70 & 41 & 60 & 65 & 55  \\
$t(3)  $   & 28 & 35 & 34 & 22 & 20 & 25 & 23 & 12  \\
$t(4)  $   & 16 & 18 & 14 &  6 & 10 & 10 &  6 &  3  \\
$t(5)  $   & 11 & 12 &  9 &  4 &  6 &  6 &  3 &  1  \\
$t(10) $   &  9 & 10 &  3 &  0 &  4 &  5 &  1 &  0  \\
\hline
\end{tabular}\label{tbl:sim-alternative-mle-H2-100-1.8}
\end{center}
\end{table} 

\begin{table}[h]
\begin{center}
\caption[percentage]{Power of $D_{n,\kappa}$ under $H_2$ ($\alpha=1.8$,\ significance
  levels~$\xi = 0.1,\ 0.05,\ n=200$)}
\ \\
\begin{tabular}{|l||cccc||cccc|}   
\hline
\multicolumn{1}{|c||}{$\xi$}  &
\multicolumn{4}{c||}{$0.1$} &
\multicolumn{4}{c|}{$0.05$} \\
 \hline
\multicolumn{1}{|c||}{$\kappa$}   & 1.0 & 2.5 & 5.0 & 10.0 & 1.0 & 2.5 & 5.0 & 10.0 \\ \hline
$N(0,2)$   & 20 & 23 & 18 &  0 & 10 & 12 &  4 &  0  \\
$t(1)  $   & 99 & 100 &100 & 100& 97 & 99 &100 &100  \\
$t(2)  $   & 83 & 94 & 94 & 89 & 72 & 88 & 89 & 80  \\
$t(3)  $   & 45 & 54 & 49 & 33 & 28 & 37 & 37 & 21  \\
$t(4)  $   & 24 & 25 & 18 &  7 & 13 & 15 &  9 &  4  \\
$t(5)  $   & 16 & 17 &  9 &  3 &  7 &  9 &  4 &  1  \\
$t(10) $   & 13 & 12 &  9 &  1 &  5 &  6 &  2 &  0  \\
\hline
\end{tabular}\label{tbl:sim-alternative-mle-H2-200-1.8}
\end{center}
\end{table} 
\end{document}